\documentstyle{amsppt}
\NoRunningHeads \magnification=\magstephalf \hcorrection{.5truein}
\NoBlackBoxes \refstyle{C} \loadbold

\font\teneuf=eufm10 
\font\seveneuf=eufm7 
\font\fiveeuf=eufm5
\newfam\euffam \def\euf{\fam\euffam\teneuf}
\textfont\euffam=\teneuf \scriptfont\euffam=\seveneuf
\scriptscriptfont\euffam=\fiveeuf %
\def\Afr{{\euf A}}

\def\Pfr{{\euf P}}       
\def\Qfr{{\euf Q}}



 \overfullrule=0pt

\def\thfill{\null\nobreak\hfill}
\def\endbeweis{\thfill\vbox{\hrule
  \hbox{\vrule\hbox to 5pt{\vbox to 5pt{\vfil}\hfil}\vrule}\hrule}}

\def\ref#1{\par\noindent\hangindent3\parindent
\hbox to 3\parindent{#1\hfil}\ignorespaces}

\define\K{\bold K}
\define\A{\bold A}
\define\F{\bold F}
\define\Pbf{\bold P}
\define\E{\bold E}

\def\Bscr{{\Cal B}}

\def\lra{\longrightarrow}

\def\Fix{\mathop{\roman{Fix}}\nolimits}
\def\Gal{\mathop{\roman{Gal}}\nolimits}
\def\Ker{\mathop{\roman{Ker}}\nolimits}
\def\Adm{\mathop{\roman{Adm}}\nolimits}
\def\Perm{\mathop{\roman{Perm}}\nolimits}

\def\ol{\overline}
\def\pn{\par\noindent}
\def\inv{\mathop{\roman{inv}}\nolimits}
\def\id{\mathop{\roman{id}}\nolimits}

\topmatter

\title A guide to the reduction modulo $p$ of Shimura varieties\endtitle

  \author
 M. Rapoport
\footnote{ Mathematisches Institut der Universit\"at zu K\"oln, Weyertal
86-90, D--50931 K\"oln, Deutschland}
  \endauthor

\endtopmatter

\define\C{\bold C}
\define\Dbf{\bold D}
\define\G{\bold G}

\define\Q{\bold Q}
\define\R{\bold R}
\define\Z{\bold Z}

\document


This report is based on my lecture at the Langlands conference in
Princeton in 1996 and the series of lectures I gave at the semestre Hecke
in Paris in 2000. In putting the notes for these lectures in order, it was
my original intention to give a survey of the activities in the study of
the reduction of Shimura varieties. However, I realized very soon that
this task was far beyond my capabilities. There are impressive results on
the reduction of ``classical'' Shimura varieties, like the Siegel spaces
or the Hilbert-Blumenthal spaces, there are deep results on the reduction
of specific Shimura varieties and their application to automorphic
representations and modular forms, and to even enumerate all these
achievements of the last few years in one report would be very difficult.
Instead, I decided to concentrate on the reduction modulo $p$ of Shimura
varieties {\it for a parahoric level structure} and more specifially on
those aspects which have a {\it group-theoretic interpretation.} Even in
this narrowed down focus it was not my aim to survey all results in this
area but rather to serve as a guide to those problems with which I am
familiar, by putting some of the existing literature in its context and by
pointing out unsolved questions. These questions or conjectures are of two
different kinds. The first kind are open even for those Shimura varieties
which are moduli spaces of abelian varieties. Surely these conjectures are
the most urgent and the most concrete and the most tractable. The second
kind are known for these special Shimura varieties. Here the purpose of
the conjectures resp.\ questions is to extend these results to more
general cases, e.g.\ to Shimura varieties of Hodge type.

As a general rule, I wish to stress that I would not be surprised if some
of the conjectures stated here turn out to be false, especially in cases
of very bad ramification. But I believe that even in these cases I should
not be too far off the mark, and that a suitable modification of these conjectures
gives the correct answer. My motivation in running the risk of stating
precise conjectures is that I wanted to point out directions of
investigation which seem promising to me.

The guiding principle of the whole theory presented here is to give a
group-theoretical interpretation of phenomena found in special cases in a
formulation which makes sense for a general Shimura variety. This is
illustrated in the first section which treats some aspects of the elliptic
modular case from the point of view taken in this paper. The rest of the
article consists of two parts, the local theory and the global theory.
Their approximate contents may be inferred from the table of contents
below.

I should point out that the development in these notes is very uneven and
that sometimes I have gone into the nitty gritty detail, whereas at other
times I only give a reference for further developments. My motivation for
this is that I wanted to give a real taste of the whole subject --- in the
hope that it is attractive enough for a student, one motivated enough to
read on and skip parts which he finds unappealing.

In conclusion, I would like to stress, as in the introduction of [R2], the
influence of the ideas of V.~Drinfeld, R.~Kottwitz, R.~Langlands and
T.~Zink on my way of thinking about the circle of problems discussed here.
In more recent times I also learned enormously from G.~Faltings,
A.~Genestier, U.~G\"ortz, J.~de~Jong, E.~Landvogt, G.~Laumon, B.C.~Ng\^o,
G.~Pappas, H.~Reimann, H.~Stamm, and T.~Wedhorn, but the influence of R.~Kottwitz continued to be all-important.
 I am happy to express my
gratitude to all of them. I also thank T.~Ito, R.~Kottwitz and especially T.~Haines for their
remarks on a preliminary version of this paper.

\bigskip

\bigskip
\centerline{\bf Table of contents}

\bigskip

\itemitem{1.} {Motivation: The elliptic modular curve}

\medskip
\item{I.} {Local theory}

\medskip
\itemitem{2.} {Parahoric subgroups}
\itemitem{3.} {$\mu$-admissible and $\mu$-permissible set}
\itemitem{4.} {Affine Deligne-Lusztig varieties}
\itemitem{5.} {The sets $X(\mu, b)_K$}
\itemitem{6.} {Relations to local models}

\bigskip
\item{II.} {Global theory}

\medskip
\itemitem{7.} {Geometry of the reduction of a Shimura variety}
\itemitem{8.} {Pseudomotivic and quasi-pseudomotivic Galois gerbs}
\itemitem{9.} {Description of the point set in the reduction}
\itemitem{10.} {The semi-simple zeta function}

\medskip
\itemitem{} {Bibliography}

\bigskip

\subheading{1. Motivation: The elliptic modular curve}

\bigskip

In this section we illustrate the problem of describing the reduction
modulo $p$ of a Shimura variety in the simplest case. Let ${\G}=GL_2$ and
let $({\G},\{ h\})$ be the usual Shimura datum. Let $\K\subset \G(\A_f)$
be an open compact subgroup of the form $\K =K^p.K_p$ where $K^p$ is a
sufficiently small open compact subgroup of ${\G}(\A_f^p)$. Let
$G=\G\otimes_{\Q}\Q_p$. We consider the cases where $K_p$ is one of the following two parahoric
subgroups of $G(\Q_p)$,

\medskip\noindent
(i) $K_p=K_p^{(i)}=GL_2(\Z_p)$ ({\it hyperspecial maximal parahoric})

\medskip\noindent
(ii) $K_p=K_p^{(ii)}=\left\{ g\in GL_2(\Z_p);\ g\equiv \pmatrix *&*\cr
0&*\cr\endpmatrix {\roman{mod}}\ p\right\}$ ({\it Iwahori})

\medskip\noindent
The corresponding Shimura variety $Sh(\G, h)_{\K}$ is defined over $\Q$.
It admits a model $Sh(G,h)_{\K}$ over ${\roman{Spec}}\ \Z_{(p)}$ by posing
the following moduli problem over $(Sch/\Z_{(p)})$:

\medskip\noindent
(i) an elliptic curve $E$ with a level-$K^p$-structure.

\medskip\noindent
(ii) an isogeny of degree $p$ of elliptic curves $E_1\to E_2$, with a
level-$K^p$-structure.

\medskip\noindent
The description of the point set $Sh(G,h)_{\K}(\overline\F_p)$ takes in
both cases (i) and (ii) the following form,
$$Sh(G,h)_{\K}(\overline\F_p)=\coprod\limits_{\varphi}
I_{\varphi}(\Q)\setminus X(\varphi)_{K_p}\times X^p/K^p\ \ .\leqno(1.1)$$
Here the sum ranges over the isogeny classes of elliptic curves and
$I_{\varphi}(\Q)= {\roman{End}}_{\Q}(E)^\times$ is the group of
self-isogenies of any element of this isogeny class. Furthermore,
$X^p/K^p$ may be identified with $\G(\A_f^p)/K^p$, with the action of
$I_{\varphi}(\Q)$ defined by the $\ell$-adic representation afforded by
the rational Tate module. The set $X(\varphi)_{K_p}$ is the most
interesting ingredient.

Let ${\Cal O}=W(\overline\F_p)$ be the ring of Witt vectors over
$\overline\F_p$ and $L={\roman{Fract}}\ {\Cal O}$ be its fraction field.
We denote by $\sigma$ the Frobenius automorphism of $L$. Let $N$ denote
the rational Dieudonn\'e module of $E$. Then $N$ is a 2-dimensional
$L$-vector space, equipped with a $\sigma$-linear bijective endomorphism
$F$ (the crystalline Frobenius). Then in case (i) (hyperspecial case), the
set $X(\varphi)_{K_p^{(i)}}$ has the following description $$\align
X(\varphi)_{K_p^{(i)}} & =\{\Lambda;\
p\Lambda\mathop{\subset}\limits_{\neq} F\Lambda
\mathop{\subset}\limits_{\neq}\Lambda\}\tag1.2
\\
& = \{\Lambda ; {\roman{inv}}(\Lambda, F\Lambda)=\mu\}\ .
\endalign$$
Here $\Lambda$ denotes a ${\Cal O}$-lattice in $N$. The set of ${\Cal
O}$-lattices in $N$ may be identified with $G(L)/G({\Cal O})$. We have
used the elementary divisor theorem to establish an identification
$${\roman{inv}}:G(L)\setminus [G(L)/G({\Cal O})\times G(L)/G({\Cal O})]
=G({\Cal O})\setminus G(L)/G({\Cal O})\simeq \Z^2/S_2\ .\leqno(1.3)$$
Furthermore $\mu=(1,0)\in\Z^2/S_2$ is the conjugacy class of one-parameter
subgroups associated to $\{ h\}$.

In case (ii) (Iwahori case), the set $X(\varphi)_{K_p^{(ii)}}$ has the
following description, $$X(\varphi)_{K_p^{(ii)}}= \{ p\Lambda_2
\mathop{\subset}\limits_{\neq} \Lambda_1 \mathop{\subset}\limits_{\neq}
\Lambda_2;\
p\Lambda_1\mathop{\subset}\limits_{\neq}F\Lambda_1\mathop{\subset}\limits_{\neq}
\Lambda_1, p\Lambda_2\mathop{\subset}\limits_{\neq} F\Lambda_2
\mathop{\subset}\limits_{\neq} \Lambda_2\}\ .\leqno(1.4)$$ Here again
$\Lambda_1,\Lambda_2$ denote ${\Cal O}$-lattices in $N$.

In either case $X(\varphi)_{K_p}$ is equipped with an operator $\Phi$
which under the bijection (1.1) corresponds to the action of the Frobenius
automorphism on the left hand side.

Let us describe the set $X(\varphi)_{K_p^{(ii)}}$ in the manner of the
second line of (1.2). The analogue in this case of the relative position
of two chains of inclusions of ${\Cal O}$-lattices in $N$, $p\Lambda_2
\mathop{\subset}\limits_{\neq}\Lambda_1 \mathop{\subset}\limits_{\neq}
\Lambda_2$ and $p\Lambda'_2\mathop{\subset}\limits_{\neq}
\Lambda'_1\mathop{\subset}\limits_{\neq} \Lambda'_2$ is given by the
identification analogous to (1.3), $${\roman{inv}}:G(L)\setminus [G(L)/
G_0({\Cal O})\times G(L) /G_0 ({\Cal O})] = G_0 ({\Cal O})\setminus
G(L)/G_0({\Cal O})\buildrel\sim\over\to\Z^2\rtimes S_2.\leqno(1.5)$$ Here
$G_0({\Cal O})$ denotes the standard Iwahori subgroup of $G({\Cal O})$ and
on the right appears the extended affine Weyl group $\tilde W$ of $GL_2$.
It is now a pleasant exercise in the Bruhat-Tits building of $PGL_2$ to
see that $$\align & \{ p\Lambda_2 \mathop{\subset}\limits_{\neq} \Lambda_1
\mathop{\subset}\limits_{\neq} \Lambda_2,\
p\Lambda'_2\mathop{\subset}\limits_{\neq}
\Lambda'_1\mathop{\subset}\limits_{\neq} \Lambda'_2; \tag1.6\\ &
p\Lambda_1\mathop{\subset}\limits_{\neq}
\Lambda'_1\mathop{\subset}\limits_{\neq} \Lambda_1,\
p\Lambda_2\mathop{\subset}\limits_{\neq}
\Lambda'_2\mathop{\subset}\limits_{\neq} \Lambda_2\}\\
=
& \{ (g,g')\in (G(L)/ G_0({\Cal O}))^2;\ {\roman{inv}} (g,g')\in
{\roman{Adm}}(\mu)\}\ .
\endalign$$
Here ${\roman{Adm}}(\mu)$ is the following subset of $\tilde W$,
$${\roman{Adm}}(\mu)=\{ t_{(1,0)}, t_{(0,1)}, s\cdot t_{(1,0)}\}\
.\leqno(1.7)$$ Here $t_{(1,0)}$ and $t_{(0,1)}$ denote the translation
elements in $\tilde W=\Z^2\rtimes S_2$ corresponding to (1,0) resp.\ (0,1)
in $\Z^2$, and $s$ denotes the non-trivial element in $S_2$.

For $w\in\tilde W$ let us introduce the {\it affine Deligne-Lusztig
variety,} $$X_w(F)= \{ g\in G(L)/G_0 ({\Cal O});\ {\roman{inv}}(g,Fg)=w\}\
.\leqno(1.8)$$ Then we may rewrite (1.4) in the following form,
$$X(\varphi)_{K_p^{(ii)}}=\bigcup\limits_{w\in{\roman{Adm}}(\mu)} X_w(F)\
\ .\leqno(1.9)$$ This is analogous to the second line in (1.2) which may
be viewed as a {\it generalized} affine Deligne-Lusztig variety
corresponding to the hyperspecial parahoric $K_p^{(i)}$. It should be
pointed out that in this special case the union (1.9) is spurious: only
one of the summands is non-empty, for a fixed isogeny class $\varphi$. For
more general Shimura varieties this is no longer true.

The model $Sh(G,h)_{\K}$ is smooth over ${\roman{Spec}}\ \Z_{(p)}$ in the
hyperspecial case, but it has bad reduction in the Iwahori case. In the
latter case there is the famous picture of the special fiber where two
hyperspecial models meet at the supersingular points.

Such a global picture is not known in more general cases. The nature of
the singularities in the special fiber in the Iwahori case can be
understood in terms of the associated {\it local model.}

We consider the lattice chain $p\Lambda_2\mathop{\subset}\limits_{\neq}
\Lambda_1\mathop{\subset}\limits_{\neq} \Lambda_2$ in $\Q_p^2$, where
$\Lambda_2= \Z_p^2$ and $\Lambda_1=p\Z_p\oplus \Z_p$. Let ${\Cal
M}^{\roman{loc}}(G,\mu)_{K_p^{(ii)}}$ be the join of $\Pbf(\Lambda_1)$ and
$\Pbf(\Lambda_2)$ over $\Z_p$ (= scheme-theoretic closure of the common
generic fiber $\Pbf^1_{\Q_p}$ in $\Pbf(\Lambda_1)\times_{\roman{Spec}\
\Z_p} \Pbf(\Lambda_2))$. Then we obtain a diagram of schemes over
${\roman{Spec}}\ \Z_p$

$$\matrix \widetilde{Sh}(G,h)_{\K}\\
\llap{$\scriptstyle\pi$}\big\downarrow & \searrow\rlap{$\scriptstyle
\tilde\lambda$}\\ Sh(G,h)_{\K} && {\Cal
M}^{\roman{loc}}(G,\mu)_{K_p^{(ii)}}\ \ .
\endmatrix
\leqno(1.10)$$ Here $\pi$ is the principal homogeneous space under the
group scheme ${\Cal G}$ over ${\roman{Spec}}\ \Z_p$ attached to
$K_p^{(ii)}$ (with ${\Cal G}(\Z_p)=K_p^{(ii)}$), which adds to the isogeny
of degree $p$, $E_1\to E_2$, and its level-$K^p$-structure, a
trivialization of the DeRham homology modules, $$\matrix H_{DR}(E_1) &
\longrightarrow & H_{DR}(E_2)
\\
\big\downarrow\rlap{$\scriptstyle \simeq$} &&
\big\downarrow\rlap{$\scriptstyle \simeq$}
\\
\Lambda_1\otimes_{\Z_p}{\Cal O}_S & \longrightarrow &
\Lambda_2\otimes_{\Z_p}{\Cal O}_S & .
\endmatrix
\leqno(1.11)$$ The morphism $\tilde\lambda$ is given by the Hodge
filtration of the DeRham homology.

The diagram (1.10) can be interpreted as a relatively representable
morphism of algebraic stacks $$\lambda: Sh(G,h)_{\K}\longrightarrow [{\Cal
M}^{\roman{loc}}(G,\mu)_{K_p^{(ii)}}/{\Cal G}]\ \ .\leqno(1.12)$$ This
morphism $\lambda$ is smooth of relative dimension ${\roman{dim}}\ {\Cal
G}=4$. The analogue of $\lambda$ in the hyperspecial case is a smooth
morphism of relative dimension 4, $$\lambda: Sh(G,h)_{\K}\longrightarrow
[{\bold P}^1_{\Z_p}/ {GL_2}_{\Z_p}]\ \ .\leqno(1.13)$$

At this point we have met in this special case all the main actors which
will appear in the sequel: the admissible subset of the extended affine
Weyl group, affine Deligne-Lusztig varieties, the sets $X(\varphi)_{K_p}$
(later denoted by $X(\mu,b)_K)$, local models etc. These definitions can
be given purely in terms of the $p$-adic group $GL_2$ and its parahoric
subgroup $K_p$. This will be the subject matter of the local part
(sections 2--6). On the other hand, the enumeration of isogeny classes and
the description (1.1) of the points in the reduction are global problems.
These are addressed in the global part (sections 7--10).

We conclude this section with the definition of a {\it Shimura variety of
PEL-type}. The guiding principle of the theory is to investigate the
moduli problems related to them and then to express these findings in
terms of the Shimura data associated to them.

Let $B$ denote a finite-dimensional semi-simple $\Q$-algebra, let $*$ be a
positive involution on $B$, let $V\neq (0)$ be a finitely generated left
$B$-module and let $\langle\ ,\ \rangle$ be a non-degenerate alternating
bilinear form $\langle\ ,\ \rangle : V\times V\to\Q$ of the underlying
$\Q$-vector space such that $\langle bv, w\rangle =\langle v, b^*w\rangle$
for all $v,w\in V$, $b\in B$. We denote by $G$ the group of $B$-linear
symplectic similitudes of $V$. This is an algebraic group $\G$ over $\Q$.
We assume that $\G$ is a {\it connected,} hence reductive, algebraic
group. We let $h:\C^\times\to \G(\R)$ be an algebraic homomorphism which
defines on $V_{\R}$ a Hodge structure of type $(-1,0)+(0,-1)$ and which
satisfies the usual Riemann conditions with respect to $\langle\ ,\
\rangle$, comp.\ [W], 1.3.1. These data define by Deligne a Shimura
variety $Sh(\G, h)$ over the Shimura field ${\bold E}$.

We now fix a prime number $p$. Let $G=\G\otimes_{\Q}\Q_p$. We consider an
order $O_B$ of $B$ such that $O_B\otimes\Z_p$ is a maximal order of
$B\otimes \Q_p$. We assume that $O_B\otimes \Z_p$ is invariant under the
involution $*$. We also fix a multichain ${\Cal L}$ of $O_B\otimes
\Z_p$-lattices in $V\otimes \Q_p$ which is self-dual for $\langle\ ,\
\rangle$, [RZ2]. Then the stabilizer of ${\Cal L}$ is a parahoric subgroup
$K=K_p$ of $G(\Q_p)$.

Finally we fix an open compact subgroup $K^p\subset \G(\A_f^p)$, which
will be assumed sufficiently small. Let ${\bold K} = K^p.K_p$.

We fix embeddings $\overline\Q\to\C$ and $\overline\Q\to\overline\Q_p$. We
denote by ${\bold p}$ the corresponding place of ${\bold E}$ over $p$ and
by $E={\bold E}_{\bold p}$ the completion and by $\kappa$ the residue
field of ${\Cal O}_E$. We then define a moduli problem $Sh(G,h)_{\K}$ over
${\roman{Spec}}\ {\Cal O}_{E}$, i.e., a set-valued functor, as follows. It
associates to a scheme $S$ over ${\Cal O}_{E}$ the following data up to
isomorphism ([RZ2], 6.9).

\medskip
\item{1)} An ${\Cal L}$-set of abelian varieties $A=\{ A_\Lambda;\ \Lambda\in
{\Cal L}\}$.

\item{2)} A $\Q$-homogeneous principal polarization $\overline\lambda$ of
the ${\Cal L}$-set $A$.

\item{3)} A $K^p$-level structure
$$\overline\eta: H_1(A, \A_f^p)\simeq V\otimes \A_f^p\ {\roman{mod}}\ K^p\
\ ,$$ which respects the bilinear forms on both sides up to a constant in
$(\A_f^p)^\times$.

\par\noindent
We require an identity of characteristic polynomials for each $\Lambda\in
{\Cal L}$, $${\roman{char}}(b;\ {\roman{Lie}}\ A_\Lambda)=
{\roman{char}}(b; V_h^{0,-1})\ \ ,\ \ b\in O_B\ \ .$$ This moduli problem
is representable by a quasi-projective scheme whose generic fiber is the
initial Shimura variety $Sh(\G, h)_{\bold K}$ (or at least a finite union
of isomorphic copies of $Sh(\G, h)_{\bold K}$).

However, contrary to the optimistic conjecture in [RZ2], this does not
always provide us with a good integral model of the Shimura variety, e.g.\
flatness may fail. However, if the center of $B$ is a product of field
extensions which are unramified at $p$, then $Sh(G,h)_{\K}$ is a good
integral model of the Shimura variety [G1], [G2]. For most of the
remaining cases there is a closed subscheme of the above moduli space
which is a good model [PR1], [PR2]. However, these closed subschemes
cannot be defined in terms of the moduli problem of abelian varieties.
Still, they can be analyzed and can serve as an experimental basis for the
predictions which are the subject of this report.

\bigskip

\heading{I.} {Local theory}
\endheading
\bigskip
\subheading{2. Parahoric subgroups}

\bigskip

Let $G$ be a connected reductive group over a complete discretely valued
field $L$ with algebraically closed residue field. Kottwitz [K4] defines a
functorial surjective homomorphism $$\tilde\kappa_G: G(L)\lra X^*(\hat
Z(G)^I)\ \ .\tag2.1$$ Here $I=\Gal(\ol L/L)$ denotes the absolute Galois
group of $L$ and $\hat Z(G)$ denotes the center of the Langlands dual
group. For instance, if $G=GL_n$, then the target group is $\Z$ and
$\tilde\kappa_G(g)={\roman{ord}}\ {\roman{det}}\ g$; if $G=GSp_{2n}$, then
again the target group is $\Z$ and $\tilde\kappa_G(g)={\roman{ord}}\
c(g)$, where $c(g)\in L^\times$ is the multiplier of the symplectic
similitude $g$.

Let $\Bscr=\Bscr(G_{ad}, L)$ denote the Bruhat-Tits building of the
adjoint group over $L$. Then $G(L)$ acts on $\Bscr$.
\medskip\noindent
{\bf Definition 2.1.} The {\it parahoric subgroup associated to a facet}\/
$F$ of $\Bscr$ is the following subgroup of $G(L)$, $$K_F=\Fix(F)\cap\Ker
\tilde\kappa_G\ \ .$$ If $F$ is a maximal facet, i.e.\ an alcove, then the
parahoric subgroup is called an {\it Iwahori subgroup.}\/

\medskip\noindent
{\bf Remarks 2.2.} (i) We have $$K_{gF}=gK_Fg^{-1}\ \ ,\ \ g\in G(L)\ \
.$$ In particular, since all alcoves are conjugate to each other, all
Iwahori subgroups are conjugate.
\par\noindent
(ii) This notion of a parahoric subgroup coincides with the one by Bruhat
and Tits [BT2], 5.2.6., cf.\ [HR]. Let ${\Cal O}_L$ be the ring of
integers in $L$. There exists a smooth group scheme ${\Cal G}_F$ over
${\roman{Spec}}\ {\Cal O}_L$, with generic fiber equal to $G$ and with
connected special fiber such that $$K_F={\Cal G}_F({\Cal O}_L)\ \ .$$
(iii) Let $G=T$ be a torus. Then there is precisely one parahoric subgroup
$K$ of $T(L)$. Then $$K={\Cal T}^0({\Cal O}_L)\ \ .$$ Here ${\Cal T}^0$
denotes the identity component of the N\'eron model of $T$.
\medskip\noindent
2.3. Let $S$ be a maximal split torus in $G$ and $T$ its centralizer.
Since by Steinberg's theorem $G$ is quasi-split, $T$ is a maximal torus.
Let $N=N(T)$ be the normalizer of $T$. Let $$\tilde\kappa_T:T(L)\lra
X^*(\hat T^I)= X_*(T)_I\tag2.2$$ be the Kottwitz homomorphism associated
to $T$ and let $T(L)_1$ be its kernel. The factor group $$\tilde
W=N(L)/T(L)_1\tag2.3$$ will be called the {\it Iwahori Weyl group
associated to} $S$. It is an extension of the relative Weyl group
$$W_0=N(L)/T(L)\ \ .\tag2.4$$ Namely, we have an exact sequence induced by
the inclusion of $T(L)_1$ in $T(L)$,  $$0\lra X_*(T)_I\lra \tilde W\lra
W_0\lra 1\ \ .\tag2.5$$ The reason for the name given to $\tilde W$ comes
from the following fact [HR].
\medskip\noindent
{\bf Proposition 2.4.} {\it Let $K_0$ be the Iwahori subgroup associated
to an alcove contained in the apartment associated to the maximal split
torus $S$. Then $$G(L)=K_0.N(L).K_0\ \ ,$$ and the map $K_0nK_0\longmapsto
n\in \tilde W$ induces a bijection  $$K_0\setminus G(L)/K_0\simeq \tilde
W\ \ .$$ More generally, let $K$ and $K'$ be parahoric subgroups
associated to facets contained in the apartment corresponding to $S$. Let
$$\tilde W^K=(N(L)\cap K)/T(L)_1,\ \hbox{resp.}\ \tilde W^{K'}=(N(L)\cap
K')/T(L)_1 \ .$$ Then } $$K\setminus G(L)/K'\simeq \tilde
W^K\setminus\tilde W /\tilde W^{K'}\ \ .$$

\medskip\noindent
For the structure of $\tilde W$ we have the following fact [HR].
\medskip\noindent
{\bf Proposition 2.5.} {\it Let $x$ be a special vertex in the apartment
corresponding to $S$, and let $K$ be the corresponding parahoric subgroup.
The subgroup $W^K$ of $\tilde W$ projects isomorphically to the factor
group $W_0$ and the exact sequence (2.5) presents $\tilde W$ as a
semidirect product, } $$\tilde W=W_0\ltimes X_*(T)_I\ \ .$$

\medskip\noindent
Sometimes for $\nu\in X_*(T)_I$ we write $t_{\nu}$ for the corresponding
element of $\tilde W$.

Let $S_{sc}$ resp.\ $T_{sc}$ resp.\ $N_{sc}$ be the inverse images of
$S\cap G_{der}$ resp.\ $T\cap G_{der}$ resp. $N\cap G_{der}$ in the simply
connected covering $G_{sc}$ of the derived group $G_{der}$. Then $S_{sc}$
is a maximal split torus of $G_{sc}$, and $T_{sc}$ resp. $N_{sc}$ is its
centralizer resp.\ normalizer. Hence
$$W_a= N_{sc}(L) / T_{sc}(L)_1\leqno(2.6)$$
is the Iwahori Weyl group of $G_{sc}$. This group is called the {\it
affine Weyl group associated to} $S$, for the following reason. Let us fix
a special vertex $x$ in the apartment corresponding to $S$. Then there
exists a reduced root system $~^x\Sigma$ such that Proposition 2.5
(applied to $G_{sc}$ instead of $G$) presents $W_a$ as the affine Weyl
group associated (in the sense of Bourbaki) to $~^x\Sigma$,
$$W_a= W(~^x\Sigma)\ltimes Q^\vee (~^x\Sigma)\ \ ,\leqno(2.7)$$
cf.\ [T], 1.7, comp.\ also [HR]. In other words, we have an identification
$W_0\simeq W(~^x\Sigma)$ and $X_*(T_{sc})_I\simeq Q^\vee (~^x\Sigma)$
compatibly with the semidirect product decompositions (2.7) and
Proposition 2.5. In particular, $W_a$ is a Coxeter group.

There is a canonical injective homomorphism $W_a\to \tilde W$ which
induces an injection from $X_*(T_{sc})_I$ into $X_*(T)_I$. In fact, $W_a$
is a normal subgroup of $\tilde W$ and $\tilde W$ is an extension,
$$1\to W_a\to \tilde W\to
X_*(T)_I/X_*(T_{sc})_I\to 1\ \ .\leqno(2.8)$$
The affine Weyl group $W_a$ acts simply transitively on the set of alcoves
in the apartment of $S$, cf.\ [T], 1.7. Since $\tilde W$ acts transitively on the
set of these alcoves and $W_a$ acts simply transitively,
$\tilde W$ is the semidirect product of $W_a$ with the normalizer $\Omega$
of a base alcove, i.e.\ the subgroup of $\tilde W$ which preserves the
alcove as a set, $$\tilde W=W_a\rtimes\Omega\ \ .\leqno(2.9)$$ In the
sequel we will often identify $\Omega$ with $X_*(T)_I/X_*(T_{sc})_I$.
\medskip\noindent
{\bf Remarks 2.6.} (i) Let $K=K_F$ be a parahoric subgroup and ${\Cal G}
={\Cal G}_F$ the corresponding group scheme over ${\roman{Spec}}\ {\Cal
O}_L$, cf.\ Remarks 2.2, (ii). Then $\tilde W^K$ can be identified with
the Weyl group of the special fiber $\overline{\Cal G}={\Cal
G}\otimes_{{\Cal O}_L}k$ of the group scheme ${\Cal G}$.
\par\noindent
(ii) Assume in Proposition 2.5 above that $x$ is a {\it hyperspecial}
vertex. In this case $S=T$ and $W_0$ is the {\it absolute} Weyl group of
$G$. In this case we have $$\tilde W=W_0\rtimes X_*(S)$$ and, since
$\tilde W^K=W_0$, $$\tilde W^K\setminus \tilde W /\tilde W^K= X_*(S)/W_0\
\ .$$

\bigskip

\subheading{3. $\mu$-admissible and $\mu$-permissible set}

\bigskip

We continue with the notation of the previous section. In particular, we
let $N$ resp.\ $T$ be the normalizer resp.\ centralizer of a maximal split
torus $S$ over $L$. Let $$W=N(\overline L) /T(\overline L)\tag3.1$$ be the
absolute Weyl group of $G$. Then $W_0=W^I$ is the set of invariants.
\par
Let $\{ \mu\}$ be a conjugacy class of cocharacters of $G$ over $\overline
L$. We denote by the same symbol the corresponding $W$-orbit in $X_*(T)$.
We associate to $\{\mu\}$ a $W_0$-orbit $\Lambda=\Lambda(\{\mu\})$ in
$X_*(T)_I$, as follows. Let $B$ be a Borel subgroup containing $T$,
defined over $L$. We denote the corresponding closed (absolute) Weyl
chamber in $X_*(T)_{\R}$ by $\overline C_B$. Let $\mu_B\in\{\mu\}$ be the
unique element in $\overline C_B$. Then the $W_0$-orbit $\Lambda$ of the
image $\lambda$ of $\mu_B$ in $X_*(T)_I$ is well-determined since any two
choices of $B$ are conjugate under an element of $W_0$.

\medskip\noindent
{\bf Lemma 3.1.} {\it All elements in $\Lambda$ are congruent modulo
$W_a$.}
\medskip\noindent
{\bf Proof:} Let us fix a special vertex $x$ and let us identify $W_0$
with $W(~^x\Sigma)$, cf.\ (2.7). We claim that for any $\lambda\in
X_*(T)_I$ and any $w\in W_0$ we have
$$w(\lambda)-\lambda\in Q^\vee (~^x\Sigma)\ \ .\leqno(3.2)$$
By induction on the length of $w$ (w.r.t.\ some ordering of the roots) we
may assume that $w=s_\alpha$ is a reflection about a simple root
$\alpha\in ~^x\Sigma$. But
$$s_\alpha(\lambda)-\lambda =-\langle \lambda, \alpha\rangle\cdot
\alpha^\vee\ \ .\leqno(3.3)$$
The assertion follows, since the image of $X_*(T)_I$ in
$X_*(T_{ad})_I\otimes \R =X_*(S_{sc})\otimes \R$ lies in the lattice of
coweights $P^\vee$ for $~^x\Sigma$ (this holds since $P^\vee$ acts simply
transitively on the set of special vertices in the apartment and these are
preserved under the subgroup $X_*(T)_I$ of $\tilde W$, comp.\
[HR]).\endbeweis

We shall denote by $\tau=\tau(\{\mu\})\in\Omega$
the common image of all elements of $\Lambda$. Let us now fix an alcove
$\bold a$ in the apartment corresponding to $S$. This defines a Bruhat
order on the affine Weyl group $W_a$ which we extend in the obvious way to
the semidirect product $\tilde W=W_a\rtimes \Omega$, cf.\ (2.9).
\par
Using this Bruhat order we can now introduce the $\mu${\it -admissible
subset of} $\tilde W$ [KR1], $${\roman{Adm}}(\mu)=\{ w\in \tilde W;\
w\leq\lambda\ \hbox{for some}\ \lambda\in\Lambda\}\ \ \tag3.4$$ Here we
consider the elements $\lambda\in X_*(T)_I$ as elements of $\tilde W$
({\it translation elements}). Note that by definition all elements in
${\roman{Adm}}(\mu)$ have image $\tau$ in $\Omega$.
\par
More generally, let ${\bold a}'$ be a facet of ${\bold a}$ and let $K$ be
the corresponding parahoric subgroup. Then the Bruhat order on $\tilde W$
induces a Bruhat order on the double coset space $\tilde W^K\setminus
\tilde W/\tilde W^K$ characterized by

$$\align & \tilde W^Kw_1\tilde W^K\leq \tilde W^Kw_2\tilde W^K
\Longleftrightarrow \tag3.5\cr & \exists w'_1\in \tilde W^K w_1\tilde W^K\
\hbox{and} \ \exists w'_2\in \tilde W^Kw_2\tilde W^K\ \hbox{such that}\
w'_1\leq w'_2\Longleftrightarrow\cr &\hbox{the same holds for $w'_1$ and
$w'_2$ the unique elements}\cr &\hbox{of minimal length in their
respective double cosets.}
\endalign$$

\par\noindent We then
define the {\it $\mu$-admissible subset of} $\tilde W^K\setminus \tilde W/\tilde
W^K$, $${\roman{Adm}}_K(\mu)=\{ w\in \tilde W^K\setminus \tilde W/\tilde
W^K;\ w\leq \tilde W^K\lambda\tilde W^K\ \ \hbox{for some}\ \lambda\in
\Lambda\}\ .\tag3.6$$ Since the element of minimal length in a double
coset is smaller than any element in it, the natural map
$${\roman{Adm}}(\mu)\longrightarrow {\roman{Adm}}_K(\mu)\tag3.7$$ is
surjective. In other words, $${\roman{Adm}}_K(\mu)=\ \hbox{image of}\
{\roman{Adm}}(\mu)\ \hbox{under}\ \tilde W\to \tilde W^K\setminus \tilde
W/\tilde W^K\ \ .\tag3.8$$
\medskip\noindent
We next introduce another subset of $\tilde W$. We first note that the
apartment in $\Bscr(G_{ad},L)$ corresponding to $S$ is a principal
homogeneous space under $X_*(S_{ad})_{\R}$. Let $\Lambda_{ad}$ be the
image of $\Lambda$ under the natural map $$X_*(T)_I\longrightarrow
X_*(T_{ad})_I\longrightarrow X_*(T_{ad})_I\otimes \R =X_*(S_{ad})_{\R} \ \
.$$ We denote by ${\Cal P}_{\mu}={\roman{Conv}}(\Lambda_{ad})$ the convex
hull of $\Lambda_{ad}$. Now we can define the $\mu${\it -permissible
subset} of $\tilde W$, $$\align {\roman{Perm}}(\mu)= & \{ w\in \tilde W;\
w\equiv\tau\ {\roman{mod}}\ W_a\ \hbox{and} \tag3.9 \cr & w(a)-a\in {\Cal
P}_{\mu},\ \hbox{for all}\ a\in{\bold a}\}\ . \cr
\endalign$$
Note that by convexity it suffices to impose the second condition in (3.9)
for the vertices $a_i$ of ${\bold a}$. Again there is a variant for a
parahoric subgroup $K$ corresponding to a facet ${\bold a}'$ of ${\bold
a}$. Since $\tilde W^K\subset W_a$ [HR], the first condition in the next
definition makes sense,
 $$\align {\roman{Perm}}_K(\mu)= & \{ w\in \tilde W^K\setminus \tilde W/\tilde W^K;\
w\equiv\tau\ {\roman{mod}}\ W_a\ \hbox{and} \tag3.10 \cr & w(a)-a\in {\Cal
P}_{\mu},\ \hbox{for all}\ a\in{\bold a}'\}\ . \cr
\endalign$$
Let us check that the second condition also depends only on the double
coset of $w$. Let us write $\tilde W=W_0\rtimes X_*(T)_I$ (corresponding
to the choice of a special vertex which defines the inclusion of $W_a$ in
$\tilde W$). If $x\in \tilde W^K$, let $x=x_0\cdot t_{\nu}$ with $x_0\in
W_0$ and $\nu\in X_*(T)_I$. Then for $a\in {\bold a}'$ we have $x(a)=a$
which implies $x_0(a)+x_0(\nu)=a$. Hence $$\align xw(a)-a &
=x_0w(a)+x_0(\nu)-a=x_0w(a)-x_0(a) \cr & =x_0(w(a)-a)\in x_0{\Cal
P}_{\mu}={\Cal P}_{\mu}\ \ ,  \cr
\endalign$$
which proves our claim.

There is a natural map $${\roman{Perm}}(\mu)\lra {\roman{Perm}}_K(\mu)\ \
.\tag3.11$$ However, in contrast to (3.8) it is not clear whether this map
is surjective.

\par
An important question is to understand the relation between the sets
${\roman{Adm}}(\mu)$ and ${\roman{Perm}}(\mu)$. In any case, the elements
$t_{\lambda}$ for $\lambda\in \Lambda$ are contained in both of them. In
fact, these elements are maximal in ${\roman{Adm}}(\mu)$. The following
fact is proved in [KR1]. We repeat the proof.

\medskip\noindent
{\bf Proposition 3.2.} {\it Let $G$ be split over $L$. Then
$$\Adm(\mu)\subset \Perm(\mu)\ \ .$$ In fact, ${\roman{Perm}}(\mu)$ is
closed under the Bruhat order.}

\medskip\noindent
Note that, by the preliminary remarks above, the second claim implies the
first. The significance of the second claim becomes more transparent when
we discuss local models in Section 6. Namely, in many cases
${\roman{Perm}}(\mu)$ is supposed to parametrize the set of Iwahori-orbits
in the special fiber of the local model and, since the latter is in these
cases a closed subvariety of an affine flag variety, it contains with an
orbit also all orbits in its closure. The question of when
${\roman{Adm}}(\mu)$ coincides with ${\roman{Perm}}(\mu)$ is closely
related to the flatness property of local models. In all cases known so
far, this flatness property was established by proving that the special
fiber is reduced and that the generic points of the irreducible components
of the special fiber are in the closure of the generic fiber. The property
${\roman{Adm}}(\mu)={\roman{Perm}}(\mu)$ is supposed to mean that the only
irreducible components of the special fiber are ``the obvious ones''
indexed by $t_{\lambda}$ for $\lambda\in \Lambda$, for which the
liftability problem should be visibly true (cf.\ G\"ortz [G1]--[G3] for
various cases).

\par
Returning to Proposition 3.2, we note that when $G$ is split over $L$, we
have $S=T$ and the action of $I$ is trivial. Furthermore $W_0=W$. The
reflections in the affine Weyl group will be written as $s_{\beta -m}$
where $\beta$ is a root in the sense of the euclidean root system and
$m\in\Z$. The proposition is a consequence of the following lemma.
\medskip\noindent
{\bf Lemma 3.3.} {\it Let ${\Cal P}$ be a $W_0$-stable convex polygon in
$X_*(S_{ad})_{\R}$. Let $x,y\in\tilde W$ with $x\leq y$. Let $v\in {\bold
a}$ and put $v_x=x(v)$, $v_y=y(v)$. Then $$\hbox{if}\ v_y\in v+{\Cal P}\ \ ,\ \
\hbox{then}\ v_x\in v+{\Cal P}\ \ .$$ }
\medskip\noindent
{\bf Proof.} We may assume that $x=s_{\beta-m}y$, with $\ell(x)<\ell(y)$.
Since $\beta-m$ separates ${\bold a}$ from $y({\bold a})$, it weakly
separates $v\in{\bold a}$ from $v_y$. Now $$\align (\beta-m)(v) &
=\beta(v)-m\tag3.12\\ (\beta-m)(v_y) & =\beta(v_y)-m\ \ .\tag3.13
\endalign$$
Hence we have 2 cases:
\par\noindent
$1^{\sevenrm st}$ case: $\beta(v)\leq m\leq \beta(v_y)$ \hfill\break
$2^{\sevenrm nd}$ case: $\beta(v_y)\leq m\leq \beta(v)$. \hfill\break Now
\par\noindent
\medskip\noindent
$$v_x=s_{\beta-m}(v_y)=v_y-[\beta(v_y)-m]\beta^{\vee}\ \ .\tag3.14$$ Hence
in either case, $v_x$ lies on the segment joining $v_y$ with
$v_y-[\beta(v_y)-\beta(v)]\cdot\beta^{\vee}$. Hence it suffices to show
that $s_{\beta}(v_y)+\beta(v)\beta^{\vee}\in v+{\Cal P}$.
\par\noindent
But $v_y=p+v,\ \hbox{with}\ p\in {\Cal P}$, hence $$\align s_{\beta}(v_y)
& =s_{\beta}(p)+s_{\beta}(v)\tag3.15\\ &
=s_{\beta}(p)+v-\beta(v)\cdot\beta^{\vee}\ \ .
\endalign$$
Hence $s_{\beta}(v_y)+\beta(v)\beta^{\vee} =s_{\beta}(p)+v\in v+{\Cal
P}$.\qed
\medskip\noindent
The converse inclusion is not true in general. We have the following
result which generalizes [KR1] valid for minuscule $\mu$.

\medskip\noindent
{\bf Theorem 3.4.} (Haines, Ngo [HN1]) {\it Let $G$ be either $GL_n$ or
$GSp_{2n}$. In the case of $GSp_{2n}$ assume that the dominant
representative of $\{\mu\}$ is a sum of minuscule dominant coweights. Then
${\roman{Adm}}(\mu)={\roman{Perm}}(\mu)$. }

\medskip\noindent
It may be conjectured that we have equality in general in Proposition 3.2,
if $\mu$ is a sum of minuscule dominant coweights. Note that, in the case
of $GL_n$, this condition on $\mu$ is automatically satisfied. On the
other hand, Haines and Ngo ([HN1]) have shown by example that for any split
group $G$ of rank $\geq 4$ not of type $A_n$, there exists a dominant
coweight $\mu$ such that ${\roman{Adm}}(\mu)\neq {\roman{Perm}}(\mu)$. In
[HN1] the result for $GSp_{2n}$ is obtained by relating the sets
${\roman{Adm}}(\mu)$ resp.\ ${\roman{Perm}}(\mu)$ with the corresponding
sets for the ``ambient'' $GL_{2n}$. It would be interesting to clarify
this relation in other cases.

In the sequel, until Proposition 3.10, we investigate the intersections of
${\roman{Adm}}(\mu)$ resp.\ ${\roman{Perm}}(\mu)$ with the translation
subgroup of $\tilde W$. These results are taken from unpublished notes of
Kottwitz, as completed by Haines. They will not be used elsewhere.

\medskip\noindent
{\bf Proposition 3.5.} (Kottwitz, Haines) {\it Let $G$ be split over $L$. Then}
$$X_*(T)\cap \Adm(\mu)= X_*(T)\cap \Perm(\mu)\ \ .$$ To prove this we need
a few more lemmas which will also be useful for other
purposes. For the time being we assume $G$ split. We denote by $R$ the set
of roots and by $R^+$ resp.\ $\Delta$ the set of positive resp.\ simple
roots for a fixed ordering.
\medskip\noindent
{\bf Lemma 3.6.} {\it Let $\nu\in X_*(S)$. Then $$\ell(t_{\nu})=\langle
\nu, 2w(\varrho)\rangle\ \ ,$$ where $\varrho ={1\over 2}\Sigma_{\alpha
>0}\alpha$ and where $w\in W$ is such that $w^{-1}(\nu)$ is dominant. In
other words, if $\nu$ is dominant then} $$\ell(t_{\nu})=\langle \nu,
2\varrho\rangle \ \ \hbox{and}\ \ \ell(t_{w(\nu)})= \ell (t_{\nu}),\
\hbox{all}\ w\in W\ \ .$$ {\bf Proof:} This is an immediate consequence of
[IM], Prop.\ 1.23.
\medskip\noindent
{\bf Lemma 3.7.} {\it For any $\beta\in R^+$ we have}
$$\ell(s_{\beta})< \langle \beta^{\vee}, 2\varrho\rangle\ \ .$$ {\bf
Proof:} It suffices to prove the weak inequality since the left hand side
is an odd integer and the right hand side an even integer (twice the height of the coroot
$\beta^{\vee}$). We use induction on $\ell(s_{\beta})$, the case
$\ell(s_{\beta})=1$ being trivial. So assume $\ell(s_{\beta})\geq 3$. We
first claim that $$\exists\ \alpha\in\Delta\ \ \hbox{such that}\ \
\ell(s_{\alpha}s_{\beta} s_{\alpha})=\ell(s_{\beta})-2\ \ .\leqno(3.16)$$
Indeed, let $\alpha\in\Delta$ such that $\ell(s_{\alpha}s_{\beta})=
\ell(s_{\beta})-1$. Then also
$\ell(s_{\beta}s_{\alpha})=\ell(s_{\beta})-1$. Hence there are 2 possible
configurations of $s_{\alpha}s_{\beta}s_{\alpha}, s_{\beta}s_{\alpha},
s_{\alpha}s_{\beta}, s_{\beta}$ in the Bruhat order.
\par\noindent
1) $s_{\alpha}s_{\beta}s_{\alpha} < (s_{\alpha}s_{\beta},
s_{\beta}s_{\alpha}) < s_{\beta}$
\par\noindent
2) $(s_{\alpha}s_{\beta}, s_{\beta}s_{\alpha}) < (s_{\beta},
s_{\alpha}s_{\beta}s_{\alpha})$
\medskip\noindent
It suffices to show that case 2) does not arise. In case 2, by Lemma 4.1
of [H3] we have $s_{\alpha}s_{\beta}s_{\alpha}=s_{\beta}$.
Hence $\beta$ and $s_{\alpha}(\beta)$ are proportional, hence
$s_{\alpha}(\beta)=\pm\beta$. The minus sign cannot occur since $\alpha$
is the only root in $R^+$ sent by $s_{\alpha}$ into $R^-$. Hence
$s_{\alpha}(\beta)=\beta$, i.e. $\langle \alpha,\beta^{\vee}\rangle
=\langle \alpha^{\vee},\beta\rangle =0$. Hence
$s_{\beta}(\alpha)=\alpha\in R^+$. But for any $\gamma\in R^+$ and $w\in
W_0$ we have $w^{-1}(\gamma)>0 \Leftrightarrow w\leq s_{\gamma}w$. Hence
$s_{\beta} < s_{\alpha}s_{\beta}$, a contradiction.
\par
Now start with $\alpha$ satisfying (3.16) and put $\beta'
=s_{\alpha}(\beta)$. Then by induction hypothesis
$$\ell(s_{\beta})-2=\ell(s_{\beta'})\leq \langle \beta^{'\vee},
2\varrho\rangle = \langle \beta^{\vee}, 2\varrho\rangle -2\langle
\beta^{\vee}, \alpha\rangle\ \ .$$ Hence it suffices to show that $\langle
\beta^{\vee}, \alpha\rangle \geq 1$. But if $\langle \beta^{\vee},
\alpha\rangle \leq 0$, then $s_{\beta}(\alpha)= \alpha -\langle
\beta^{\vee}, \alpha\rangle \beta\in R^+$. Hence, arguing as before,
$s_{\beta} < s_{\alpha}s_{\beta}$, a contradiction.\endbeweis

\medskip\noindent
{\bf Lemma 3.8.} {\it Let $\nu\in X_*(S)$ be dominant. Let $\beta\in R^+$
such that $\nu-\beta^{\vee}$ is dominant. Then $t_{\nu-\beta^{\vee}}\leq
t_{\nu}$ in the Bruhat order on $\tilde W$.} \pn [Here the Bruhat order on
$\tilde W$ is defined by the alcove ${\bold a}$ in $X_*(S_{ad})_{\R}$ with
apex 0 and bounded by hyperplanes $\alpha=0$ for $\alpha\in \Delta$.] \pn
{\bf Proof:} We use the identity $$s_{\beta -1}\cdot
s_{\beta}=t_{\beta^{\vee}}\ \ ,\ \ \beta\in R\ \ .\leqno(3.17)$$ Indeed,
this follows from the expression $$s_{\beta +k}(x)= x-\langle \beta,
x\rangle \beta^{\vee} -k\beta^{\vee}\ \ ,\ \ x\in X_*(S_{ad})_{\R}\ \
.\leqno(3.18)$$ This last identity also shows $$t_{\nu}\cdot
s_{\beta}=s_{\beta -m}\cdot t_{\nu}\ \ ,\ \ \beta\in R\ \ ,\ \ \nu\in
X_*(S)\ \ .\leqno(3.19)$$ Here $m=\langle \beta, \nu\rangle$. \pn The
assertion of the lemma follows from the following two statements. $$\align
t_{\nu}\cdot s_{\beta} & \leq t_{\nu} \tag3.20\\ t_{\nu-\beta^{\vee}} &
\leq t_{\nu}\cdot s_{\beta}\ \ . \tag3.21
\endalign$$
Let us prove (3.20), i.e. $$s_{\beta -m}\cdot t_{\nu}\leq t_{\nu}\ \ ,\ \
m=\langle \beta, \nu\rangle\ \ .\leqno(3.22)$$ It is enough to show that
$\beta-m$ separates ${\bold a}$ from $t_{\nu}({\bold a})$. But $$\align
(\beta-m)({\bold a}) & =\beta ({\bold a})-m\tag3.23\\ (\beta
-m)(t_{\nu}({\bold a})) & =\beta({\bold a})\subset [0,1]\ \ .
\endalign$$
Hence it suffices to know that $m\geq 1$. But since $\nu-\beta^{\vee}$ is
dominant we have $$\langle\beta, \nu-\beta^{\vee}\rangle \ge 0\ ,\
\hbox{i.e.}\ m=\langle \beta, \nu\rangle \geq \langle\beta,
\beta^{\vee}\rangle =2\ \ .\leqno(3.24)$$ Now let us prove (3.21).  It
follows with (3.17) that both sides of (3.21) differ by a reflection,
since $$t_{\nu-\beta^{\vee}}=t_\nu\cdot t_{\beta^\vee}^{-1} =(t_\nu\cdot
s_{\beta})\cdot s_{\beta -1}\ \ .\leqno(3.25)$$ Hence it suffices to prove
that $\ell(t_{\nu-\beta^{\vee}})
< \ell(t_{\nu}\cdot s_{\beta})$. But by Lemma 3.6 we have, since $\nu$ and
$\nu-\beta^{\vee}$ are dominant, $$\ell(t_{\nu-\beta^{\vee}})=
\ell(t_{\nu})- \langle 2\varrho, \beta^{\vee}\rangle < \ell (t_{\nu})-
\ell(s_{\beta})\leq \ell (t_{\nu}s_{\beta})\ .\leqno(3.26)$$ For the first
inequality we used Lemma 3.7. \endbeweis

\medskip\noindent
{\bf Remark 3.9} (Haines): In the course of the proof of Lemma 3.8 we
proved the chain of inequalities $$t_{\nu-\beta^{\vee}}\leq
t_{\nu}s_{\beta}\leq t_{\nu}$$ which is stronger than the assertion of the
Lemma. Here is a simpler argument for the assertion of Lemma 3.8 which
does not make use of Lemma 3.7. From (3.17) we have
$$t_{\nu-\beta^{\vee}}\cdot s_{\beta-1}s_{\beta}=t_{\nu}\ \
,\leqno(3.27)$$ hence $$t_{\nu-\beta^{\vee}}\equiv t_{\nu}s_{\beta -1}\ \
\hbox{in}\ W_0\setminus \tilde W/W_0\ \ .\leqno(3.28)$$ But from (3.19) we
have $$t_{\nu}s_{\beta -1}= s_{\beta -(m+1)}\cdot t_{\nu}\ \ .$$ But
$$\align & (\beta -(m+1))({\bold a})= \beta ({\bold a})- (m+1) \subset
(-\infty, -2]\\ & (\beta -(m+1))(t_{\nu}({\bold a})) =\beta({\bold
a})-1\subset [-1,0]
\endalign$$
Hence the same argument that was used to prove (3.20) also shows
that\hfill\break $s_{\beta-(m+1)}\cdot t_{\nu}\leq t_{\nu}$. We conclude
that $$W_0t_{\nu-\beta^{\vee}}W_0\leq W_0t_{\nu}W_0\ \ ,\leqno(3.29)$$ or
in terms of the longest elements in the respective double cosets
$$w_0t_{\nu-\beta^{\vee}}\leq w_0t_{\nu}\ \ .$$ But, quite generally, if
$\lambda$ and $\mu$ in $X_*(S)$ are dominant and $w\in W_0$ with
$t_{\lambda}\leq wt_{\mu}$, then $t_{\lambda}\leq t_{\mu}$. We prove this
by descending induction on $w$. Assume therefore that $t_{\lambda}\leq
wt_{\mu}$ and let $s$ be a simple reflection with $sw <w$. We wish to show
that $t_{\lambda}\leq swt_{\mu}$. But if this does not hold, then
$st_{\lambda}\leq swt_{\mu}$. Since $\lambda$ is dominant we have
$\ell(t_{\lambda})\leq \ell(st_{\lambda})$ since $\ell(t_{\lambda})=
\langle\lambda, 2\varrho\rangle$ (Lemma 3.6) and
$$\ell(st_{\lambda})=\sum\limits_{\alpha >0\atop s(\alpha)<0} \vert
\langle \alpha, \lambda\rangle +1\vert + \sum\limits_{\alpha
>0\atop s(\alpha)>0} \vert\langle \alpha, \lambda\rangle \vert
 \leqno(3.30)$$ ([IM], Prop.\ 1.23). Hence $t_{\lambda}\leq
 st_{\lambda}$ and therefore also $t_{\lambda}\leq swt_{\mu}$, a
 contradiction.\endbeweis

\medskip\noindent
{\bf Proof of Proposition 3.5:} Let $\nu\in X_*(T)\cap \Perm(\mu)$ and let
us prove that $t_{\nu}$ is $\mu$-admissible. Since $\nu$ is
$\mu$-permissible we have $\mu-\nu\in X_*(S_{sc})$. Let us first assume
that $\nu$ is dominant. Then $\nu\mathop{\leq}\limits^! \mu$, i.e.
$\mu-\nu$ is a non-negative sum with integer coefficients of simple
coroots. By the lemma of Stembridge, comp.\  [R3] there exists a sequence
of dominant elements $\nu_0=\nu\mathop{\leq}\limits^!
\nu_1\mathop{\leq}\limits^!\ldots \mathop{\leq}\limits^! \nu_r=\mu$, such
that each successive difference is a positive coroot. Applying Lemma 3.8
we conclude $$t_{\nu}\leq t_{\nu_1}\leq \ldots\leq t_{\mu}\ \
,\leqno(3.31)$$ hence $t_{\nu}$ is $\mu$-admissible.
\par
If $\nu$ is arbitrary there exists a conjugate under $w\in W$ which is
dominant, and any such conjugate by Lemma 3.6 has the same length. By a
general lemma of Haines [H3], Lemma 4.5, elements of $\tilde W$
which are conjugate under a simple reflection and of the same length are
simultaneously $\mu$-admissible. The result follows by induction by
writing $w$ as a product of simple reflections and using Lemma 3.6
repeatedly.\endbeweis

\medskip\noindent
In the preceding considerations we looked at the situation for an Iwahori
subgroup. Let us now make some comments on the subsets
${\roman{Adm}}_K(\mu)$ and ${\roman{Perm}}_K(\mu)$ of\hfill\break $\tilde
W^K\setminus \tilde W /\tilde W^K$ for an arbitrary parahoric subgroup
$K$, cf.\ (3.6) and (3.10). First of all we note that as a consequence of
Proposition 3.2 and the surjectivity of (3.8) we have the following
statement.

\medskip\noindent {\bf Proposition 3.10.}
{\it Let $G$ be split over $L$. Then $$\Adm_K(\mu)\subset \Perm_K(\mu)\ \
.$$ If $\Adm(\mu) =\Perm(\mu)$, and $\Perm(\mu)\to\Perm_K(\mu)$ is
surjective, then also $\Adm_K(\mu)\break =\Perm_K(\mu)$. }

\medskip\noindent
We note that, as proved in [KR1], all these statements hold true if $G$ is
equal to $GL_n$ or to $GSp_{2n}$ and $\mu$ is minuscule.

\medskip\noindent
Let now $G$ be split over $L$, and let $K$ be a special maximal parahoric
subgroup. We may take the vertex fixed by $K$ to be the origin in the
apartment. This identifies $$\tilde W^K\setminus\tilde W /\tilde W^K =
X_*(S) / W=X_*(S)\cap \ol C\ \ ,\leqno(3.32)$$ where $X_*(S)\cap \ol C$
are the dominant elements for some ordering of the roots. Let us choose
$\mu\in W(\mu)$ dominant and introduce the partial order as before (3.31)
(i.e., the difference is a sum of positive coroots).
\medskip\noindent
{\bf Proposition 3.11.} {\it Let $G$ be split over $L$ and let $K$ be a
special maximal parahoric subgroup. With the notations introduced we have}
$$\Adm_K(\mu)=\Perm_K(\mu)=\{ \nu\in X_*(S)\cap\ol C;\ \nu\leq\mu\}\ \ .$$
{\bf Proof:} Let $\nu\in\Perm_K(\mu)$. Then $\nu$ and $\mu$ have the same
image in $X_*(S)/X_*(S_{sc})$. Hence the condition on $\nu$ to be
$\mu$-permissible, which says that $t_{\nu}(0)\in {\Cal P}_{\mu}$, is
equivalent to $$\nu\mathop{\leq}\limits^!\mu\ \ .\leqno(3.33)$$
 Hence it suffices to show that (3.33) implies $t_{\nu}\leq
t_{\mu}$. But this is shown by the proof of Proposition 3.5
above.\endbeweis
\medskip\noindent
{\bf Corollary 3.12.} {\it In the situation of the previous proposition
assume that $\mu$ is minuscule. Then $\Adm_K(\mu)$ consists of one
element, namely $\mu\in X_*(T)/W$. }
\medskip\noindent
{\bf Proof:} Indeed in this case the set appearing in the statement of the
Proposition consists of $\mu$ only, cf.\ [K1].

\bigskip

\subheading{4. Affine Deligne-Lusztig varieties}

\bigskip

In this section we change notations. Let $F$ be a finite extension of
$\Q_p$ and let $L$ be the completion of the maximal unramified extension
of $F$ in a fixed algebraic closure $\overline F$ of $F$. We denote by
$\sigma$ the relative Frobenius automorphism of $L/F$. Let $G$ be a
connected reductive group over $F$ and let $\tilde G$ be the group over
$L$ obtained by base change. Let $\Bscr=\Bscr(G_{ad}, L)$ be the
Bruhat-Tits building of $\tilde G_{ad}$. The Bruhat-Tits building of
$G_{ad}$ over $F$ can be identified with the set of
$\langle\sigma\rangle$-invariants in $\Bscr$.

\par
We fix a maximal split torus $\tilde S$ of $\tilde G$ which is defined
over $F$ (such tori exist by [BT2], 5.1.12.) We also fix a facet ${\bold
a}'$ in the apartment corresponding to $\tilde S$ which is invariant under
$\langle\sigma\rangle$. Let $\tilde K=\tilde K({\bold a}')$ be the
corresponding parahoric subgroup of $\tilde G(L)$. The subgroup $K=\tilde
K\cap G(F)$ is called {\it the parahoric subgroup of} $G(F)$ {\it
corresponding to} ${\bold a}'$. [The subgroup $K$ determines ${\bold a}'$
uniquely, and hence we obtain a bijection between the set of parahoric
subgroups of $G(F)$, the set of $\sigma$-invariant parahoric subgroups of
$G(L)$ and the set of $\sigma$-invariant facets of $\Bscr$.] From Prop.\
2.4 we have a map (with obvious notation),

$${\roman{inv}}:\ \tilde G(L)/\tilde K\times \tilde G(L)/\tilde
K\lra\tilde W^{\tilde K}\setminus \tilde W/\tilde W^{\tilde K}\ \
,\tag4.1$$ in which the target space can be identified with the quotient
of the source by the diagonal action of $\tilde G(L)$.

\medskip\noindent
{\bf Definition 4.1.} Let $w\in \tilde W^K\setminus \tilde W/\tilde W^K$
and $b\in G(L)$. The {\it generalized affine Deligne-Lusztig variety
associated to} $w$ {\it and} $b$ is the set

$$X_w(b)=\{ g\in \tilde G(L)/\tilde K;\ {\roman{inv}}(b\sigma(g), g)= w\}\
\ .$$ When $\tilde K$ is an Iwahori subgroup, in which case $\tilde
W^{\tilde K}$ is trivial, i.e.\ $w\in \tilde W$, this set is called {\it
the affine Deligne-Lusztig variety associated to} $w$ {\it and} $b$.

\par\noindent
Let $$J_b(F)= \{ h\in G(L);\ h^{-1}b\sigma(h)=b\}\ \ .$$ Then $J_b(F)$
acts on $X_w(b)$ via $g\mapsto hg$.

\medskip\noindent
{\bf Remarks 4.2.} (i) If $b'\in G(L)$ is $\sigma$-conjugate to $b$, i.e.\
$b'=h^{-1}b\sigma(h)$, then the map $g\mapsto g'=hg$ induces a bijection

$$X_w(b)\buildrel\sim\over \lra X_w(b')\ \ .$$ Sometimes it is useful to
indicate the parahoric subgroup $K$ in the notation. If $K'\subset K$,
then, denoting by $w$ a representative of $w$ in $\tilde W^{\tilde
K'}\setminus \tilde W /\tilde W^{\tilde K'}$, there is a natural map
$$X_w(b)_{K'}\longrightarrow X_w(b)_K\ \ .$$ (ii) One could hope to equip
$X_w(b)$ with the structure of an algebraic variety locally of finite type
over the residue field $\F$ of ${\Cal O}_L$.

\par\noindent
(iii) The name given to this set derives from the analogue where $F$ is
replaced by the finite field $\F_q, L$ by the algebraic closure $\F$ of
$\F_q$ and where $\tilde K=B(\F)$ for a Borel subgroup of $G$ defined over
$\F_q$. Then $\tilde W$ is the geometric Weyl group of $G$ and with $b=1$
we obtain the usual Deligne-Lusztig variety associated to $w$ [DL]. In
this analogy, $J_b(F)$ becomes the group of rational points $G(\F_q)$. If
instead of a Borel subgroup we consider a conjugacy class of parabolic
subgroups defined over $\F_q$, we obtain the {\it generalized}
Deligne-Lusztig varieties, comp.\ [DM].

\medskip\noindent For the classical (generalized) Deligne-Lusztig varieties, there is a
simple formula for their dimensions. For affine Deligne-Lusztig varieties
such a formula is unknown. In fact, it is an open problem to determine the
pairs $(w,b)$ for which $X_w(b)\neq \emptyset$.

\medskip\noindent
{\bf Example 4.3.} Let $G=GL_2$ and $K=K_0=$ standard Iwahori subgroup. We
associate to $b\sigma$ its slope vector $\lambda=(\lambda_1, \lambda_2)\in
\Q^2$. Then $\lambda_1\geq \lambda_2$ with $\lambda_1+\lambda_2\in\Z$ and
$\lambda_i \in\Z$ if $\lambda_1\neq \lambda_2$. If $X_w(b)\neq \emptyset$
then the image of $w$ in

$$X_*(T)_I/X_*(S_{sc})=\Z$$ coincides with $\lambda_1+\lambda_2$. Conversely,
let us assume this and let us enumerate the $w\in \tilde W$ for which
$X_w(b)\neq\emptyset$. We distinguish cases.

\par\noindent
a) $b$ basic, i.e.\ $\lambda_1=\lambda_2$.
\par\noindent
a1) $\lambda_1+\lambda_2$ odd.
\par\noindent
In this case $X_w(b)\neq\emptyset \Leftrightarrow \ell(w)$ is even.
\par\noindent
a2) $\lambda_1+\lambda_2$ even.
\par\noindent
In this case $X_w(b)\neq \emptyset \Leftrightarrow$ either the projection
of $w$ to $W_a$ is trivial, or $\ell(w)$ is odd.
\par\noindent
b) $b$ hyperbolic, i.e.\ $\lambda_1\neq \lambda_2$.
\par\noindent
In this case $X_w(b)\neq \emptyset\Leftrightarrow$ either
$w=t_{\lambda_1-\lambda_2}$ or $\ell(w) >\ell(t_{\lambda_1-\lambda_2})$
and $\ell(w)\equiv \lambda_1-\lambda_2+1$ mod 2.

Furthermore, there is a simple formula for the dimension of $X_w(b)$.

\medskip\noindent
Before going on, we recall some definitions of Kottwitz [K2], [K4]. Let
$B(G)$ be the set of $\sigma$-conjugacy classes of elements of $G(L)$. The
homomorphism $\tilde\kappa_{\tilde G}$, cf.\ (2.1), induces a map
$$\kappa_G:B(G)\longrightarrow X^*(\hat Z(G)^{\Gamma})\ \ .\leqno(4.2)$$
Here $\Gamma= {\roman{Gal}}(\overline F/F)$ denotes the absolute Galois
group of $F$. We also recall the Newton map, $$\overline\nu
:B(G)\longrightarrow\Afr^+\ \ .\leqno(4.3)$$ Here the notation is as
follows. Let $G^*$ be the quasisplit inner form of $G$. Let $B^*$ be a
Borel subgroup of $G^*$ defined over $F$ and let $T^*$ be a maximal torus
in $B^*$. Then $\Afr=X_*(T^*)_{\R}^{\Gamma}$ and $\Afr^+$ denotes the
intersection of $\Afr$ with the positive Weyl chamber in $X_*(T^*)_{\R}$
corresponding to $B^*$. For instance, if $G=GL_n$, then the Newton map
associates to $b\in G(L)$, the slopes in decreasing order of the
isocrystal $(L^n, b\sigma)$. An element $b\in B(G)$ is called {\it basic}
if $\overline\nu_b$ is central, i.e.\ if $\overline\nu_b\in X_*(Z)_\R$.
This is the analogue for general $G$ of an isoclinic isocrystal. At the
opposite extreme of the basic elements of $B(G)$ are the unramified
elements. Namely, let $G=G^*$ be the quasisplit, and let $A$ be a maximal
split torus contained in $T^*$. Let $b\in A(L)$. Then $\overline\nu_b$ is
the unique dominant element in the conjugacy class of ${\roman{ord}}(b)\in
X_*(A)\subset \Afr$ (this follows from the functoriality of the Newton
map).

Let now $\{\mu\}$ be a conjugacy class of one-parameter subgroups of $G$.
Then $\{\mu\}$ determines a well-defined element $\mu^*$ in
$X_*(T^*)_{\R}$ lying in the positive Weyl chamber (use an inner
isomorphism of $G$ with $G^*$ over $\overline F$). Let
$$\overline\mu^*=[\Gamma:\Gamma_{\mu^*}]^{-1}\cdot\sum\limits_{\tau\in\Gamma
/\Gamma_{\mu^*}}\tau (\mu^*)\ \ .\leqno(4.4)$$ Then
$\overline\mu^*\in\Afr^+$. On the other hand, $\{\mu\}$ determines a
well-defined element $\mu^{\natural}$ of $X^*(\hat Z(G)^{\Gamma})$.

We define a finite subset $B(G,\mu)$ of $B(G)$ as the set of $b\in B(G)$
satisfying the following two conditions,

$$\align \kappa_G(b) & =\mu^{\natural}\tag4.5
\\
\overline\nu_b & \leq \overline\mu^*\ \ ,\tag4.6
\endalign$$
cf.\ [K4], section 6.
Here in (4.6) there appears the usual partial order on $\Afr^+$, for which
$\nu\leq\nu'$ if $\nu'-\nu$ is a nonnegative linear combination of simple
relative coroots.

The motivation for the definition of $B(G,\mu)$ comes from the following
fact. We return to the notation of the beginning of this section. Let us
assume that $G$ is quasisplit over $F$ and $\tilde G$ split over $L$,
i.e., $G$ is unramified. Let $K$ be a hyperspecial maximal parahoric
subgroup. Then $T=\tilde S$ and $\tilde W^{\tilde K}\setminus \tilde W/\tilde W^{\tilde K}$
can be identified with $X_*(\tilde S)/W_0$.

\medskip\noindent
{\bf Proposition 4.4.} ([RR]) {\it Let $\mu\in X_*(\tilde S)/W_0$. For
$b\in G(L)$ let $[b]\in B(G)$ be its $\sigma$-conjugacy class. Then}
$$X_{\mu}(b)\neq\emptyset\Longrightarrow [b]\in B(G,\mu)\ \ .$$ This is
the group theoretic version of Mazur's inequality between the Hodge
polygon of an $F$-crystal and the Newton polygon of its underlying
$F$-isocrystal.

\medskip\noindent
{\bf Example 4.5.} Let $G=GL_n$ and let $T=S$ be the group of diagonal
matrices, and $K$ the stabilizer of the standard lattice $O_F^n$ in $F^n$.
Then, with the choice of the upper triangular matrices for $B^*$,
$$\Afr^+=(\R^n)_+=\{ \nu= (\nu_1,\ldots, \nu_n)\in\R^n;\ \nu_1\geq
\nu_2\geq\ldots\geq \nu_n\}\ \ .$$ For the usual partial order on $\Afr^+$
we have $\nu\leq\nu'$ iff $$\sum\limits_{i=1}^r\nu_i\leq
\sum\limits_{i=1}^r\nu'_i\ \  \hbox{for}\ r=1,\ldots, n-1\ \  \hbox{and}\
\sum\limits_{i=1}^n \nu_i=\sum\limits_{i=1}^n \nu'_i\ \ .$$ Let $b\in
G(L)$ and let $(N,\Phi)=(L^n, b\cdot\sigma)$ be the corresponding
isocrystal of dimension $n$. If $M$ is an ${\Cal O}_L$-lattice in $N$ we
have $$\mu(M)={\roman{inv}} (M,\Phi(M))\in (\Z^n)_+\ \ .$$ Here
$(\Z^n)_+=\Z^n\cap (\R^n)_+$ and $\mu(M)= (\mu_1,\ldots, \mu_n)$ iff there
exists a ${\Cal O}_L$-basis $e_1,\ldots, e_n$ of $M$ such that
$\pi^{\mu_1}e_1,\ldots, \pi^{\mu_n}e_n$ is a ${\Cal O}_L$-basis of
$\Phi(M)$. By $\pi$ we denoted a uniformizer of $F$. Denoting by
$\overline\nu_b$ the Newton vector of the isocrystal $(N,\Phi)$, Mazur's
inequality states that $$\overline\nu_b\leq\mu(M)\ \ ,$$ i.e. $[b]\in
B(G,\mu(M))$.

\medskip\noindent
{\bf Conjecture 4.6.} {\it The converse in the implication of Proposition
4.4 holds.}

\medskip\noindent
In this direction we have the following results.

\medskip\noindent
{\bf Theorem 4.7.} {\it The converse implication in Proposition 4.4 holds
in either of the following cases.
\item{(i)} {\rm{[KR2]}} $G=GL_n$ or $G=GSp_{2n}$.
\item{(ii)} {\rm{[R3]}} The derived group of $G$ is simply connected, $\mu\in
X_*(\tilde S)/W_0$ is $\langle\sigma\rangle$-invariant and $b\in A(F)$,
where $A$ denotes a maximal $F$-split torus in $G$.}

\medskip\noindent
Whereas the individual affine Deligne-Lusztig varieties are very difficult
to understand, the situation seems to change radically when we form a
suitable finite union of them. This is the subject of the next section.

\bigskip

\subheading{5. The sets $X(\mu, b)_K$}

\bigskip

We continue with the notation of the previous section. Let $\{\mu\}$ be a
conjugacy class of one-parameter subgroups of $G$. Equivalently, $\{\mu\}$
is a $W$-orbit in $X_*(T)$ where $T$ denotes the centralizer of $\tilde
S$. We again introduce the subsets $\Adm_{\tilde K}(\mu)$ resp.\
$\Perm_{\tilde K}(\mu)$ of $\tilde W^{\tilde K}\setminus \tilde W /\tilde
W^{\tilde K}$.

\par
Let $b\in G(L)$. Then we define the following set, a finite union of
generalized affine Deligne-Lusztig varieties,  $$X(\mu, b)_K= \{ g\in
\tilde G(L)/\tilde K;\ \inv (g,b\sigma(g)) \in \Adm_{\tilde K}(\mu)\}\ \
.\tag5.1$$

\par
Let $E\subset\ol F$ be the field of definition of $\{\mu\}$. Let
$E_0=E\cap L$ and $r=[E_0:F]$. We note that $\sigma$ acts in compatible
way on $\tilde W$ and its subgroup $X_*(T)_I$, and that the map (4.1) is
compatible with this action.

\medskip\noindent
{\bf Lemma 5.1.} {\it The subset $\Lambda(\{\mu\})$ of $X_*(T)_I$ is
invariant under $\sigma^r$. Hence also the subsets $\Perm_{\tilde K}(\mu)$
and $\Adm_{\tilde K}(\mu)$ of $\tilde W^{\tilde K}\setminus \tilde W/\tilde W^{\tilde K}$
are invariant under $\sigma^r$. }

\medskip\noindent
{\bf Proof.} Let $\nu\in X_*(T)$ with image $[\nu]_I$ in $X_*(T)_I$. Then

$$\sigma^r([\nu]_I)= [\tau(\nu)]_I\ \ ,$$ for $\tau\in \Gal(\overline
L/F)$ an arbitrary lifting of $\sigma^r$. We take for $\tau$ an extension
of the automorphism $\id\otimes \sigma^r$ of $EL=E\otimes_{E_0}L$. Then
$\tau\in \Gal(\overline L/E)$ and hence preserves the orbit $\{\mu\}$ in
$X_*(T)$. Furthermore $\tau$ normalizes $\Gal(\overline L/L)$. Hence if
$\mu\in\overline C_B\cap\{\mu\}$ for a Borel subgroup defined over $L$,
then $\tau(\mu)\in\overline C_{B'}\cap \{\mu\}$ for another Borel subgroup
defined over $L$. Hence if $\lambda\in X_*(T)_I$ denotes the image of
$\mu$, then $\sigma^r(\lambda)=w_0(\lambda)$ for some $w_0\in W_0$ which
implies the first assertion. The second assertion follows (for
$\Adm_{\tilde K}(\mu)$ use that $w_1\leq w_2$ implies $\sigma(w_1)\leq
\sigma(w_2)$). \endbeweis

\medskip\noindent
Using this lemma we can now define an operator $\Phi$ on $X(\mu, b)_K$ by
$$\Phi(g)= (b\sigma)^r\cdot g\cdot\sigma^{-r} =b\cdot\sigma(b)\ldots
\sigma^{r-1}(b)\cdot\sigma^r(g)\ \ .\tag5.2$$ Let us check that $\Phi$
indeed preserves the set $X(\mu, b)_K$. We have $$\align \inv(\Phi(g),
b\sigma (\Phi(g)) & =\inv (\sigma^r(g), \sigma^r(b\sigma(g)))\tag5.3\\ &
=\sigma^r(\inv(g, b\sigma(g))\ \ .
\endalign$$
The claim follows from Lemma 5.1.


\medskip\noindent
In the context of Remark 4.2., (ii), the set $X(\mu, b)_K$ may be expected
to be the set of $\F$-points of an algebraic variety over $\F$, and $\Phi$
would define a Weil descent datum over the residue field $\kappa_E$ of $E$
in the sense of [RZ2].

\par
As mentioned at the end of the last section, whereas it seems difficult to
understand when the individual affine Deligne-Lusztig varieties which make
up $X(\mu, b)_K$ are non-empty, their union seems to behave better, at
least in the cases when $\{\mu\}$ is minuscule.

\medskip\noindent
{\bf Conjecture 5.2.} {\it Let $\{\mu\}$ be minuscule.
\par\noindent
a) $X(\mu, b)_K\neq\emptyset$ if and only if the class $[b]$ of $b$ in
$B(G)$ lies in the subset $B(G,\mu)$.

\par\noindent
b) For $K\subset K'$, the induced map $X(\mu, b)_K\to X(\mu, b)_{K'}$ is
surjective.}

\medskip\noindent
{\bf Remark 5.3.} (i) Suppose that $K$ is hyperspecial. Then, if $X(\mu,
b)_K\neq \emptyset$, it follows that $[b]\in B(G,\mu)$, cf.\ Prop.\ 4.4.
This holds even when $\{\mu\}$ is not minuscule. In general, it is not
clear whether the hypothesis that $\{\mu\}$ be minuscule is indeed
necessary in Conjecture 5.2.

\medskip\noindent
In the direction of Conjecture 5.2 we first note the following easy
observation.

\medskip\noindent
{\bf Lemma 5.4.} {\it If $X(\mu,
 b)_K\neq\emptyset$, then} $$\kappa(b)= \mu^{\natural}\ \ .$$

\par\noindent
{\bf Proof:} We consider the composition $\tilde{\kappa}=\tilde \kappa_G$ of
$\tilde\kappa_{\tilde G}$ and the natural surjection, $$G(L)\to X^*(\hat
Z(G)^I)\to X^*(\hat Z(G)^{\Gamma})\ \ .\leqno(5.4)$$ The map
$\tilde\kappa_G$ induces $\kappa_G$ on $B(G)$. If $g\tilde K\in X(\mu,
b)_K$, then $g^{-1} b\sigma(g)= k_1wk_1$ with $k_1, k_2\in\tilde K$ and
with $w\in {\roman{Adm}}_K(\mu)$. Since $k_1, k_2\in {\roman{Ker}}\,
\tilde\kappa_{\tilde G}$ we conclude that $$\tilde\kappa(b)=
\tilde\kappa(g^{-1}b\sigma(g)) =\tilde\kappa(k_1wk_2)= \tilde\kappa(w)\ \
.$$ But $$\tilde W^{\tilde K}w\tilde W^{\tilde K} \leq \tilde W^{\tilde K}
t_{\mu'} \tilde W^{\tilde K}\ \ ,$$ for a conjugate $\mu'$ of $\mu$.
Since $\tilde W^{\tilde K}\subset W_a$ and $W_a\subset {\roman{Ker}}\
\tilde\kappa$ (since $\tilde\kappa(G_{sc}(L)=0$), we conclude that
$\tilde\kappa(w)= \tilde\kappa(t_{\mu'})$. If $\mu'=w_0(\mu)$ for $w_0\in
W_0$ we have
$$\tilde\kappa(t_{\mu'})= \tilde\kappa(w_0t_{\mu}w_0^{-1})=
\tilde\kappa(t_{\mu})\ \ ,$$ hence
$\tilde\kappa(b)=\tilde\kappa(t_{\mu})=\mu^{\natural}$.\endbeweis

\medskip\noindent
{\bf Theorem 5.5.} ([KR2]) {\it Conjecture 5.2 holds in the cases
$G=R_{F'/F}(GL_n)$ and $G=R_{F'/F}(GSp_{2n})$, where $F'$ is an unramified
extension of $F$.}

In fact, in loc.cit.\ also the case when $G$ is an inner form
of $GL_n$ is treated.

We also mention the following case when Conjecture 5.2 holds.

\medskip\noindent
{\bf Proposition 5.6.} {\it Assume that $G$ splits over $L$ and that the center of $G$ is connected.
Let $b\in G(L)$ be such that $[b]\in B(G)_{\roman{basic}}$. Let $K_0$ be an
Iwahori subgroup defined over $F$. Then $X(\mu,
b)_{K_0}\neq\emptyset\Leftrightarrow [b]\in B(G, \mu)$. }

\medskip\noindent
{\bf Proof:} It is obvious that if $[b]\in B(G)_{basic}$, then $[b]\in
B(G,\mu)$ iff $\kappa(b)=\mu^{\natural}$. Hence one implication
$(\Rightarrow)$ follows from Lemma 5.4. Now let $[b]\in B(G,\mu)\cap
B(G)_{basic}$.

\medskip\noindent
{\it Claim:} {\it Let $N$ be the normalizer of a maximal torus $S$ which
splits over $L$, and let $\tilde K\subset G(L)$ be a $\sigma$-stable
Iwahori subgroup corresponding to an alcove in the apartment of $S$. Then
there exists a representative $b'$ of $[b]$ in $N(L)$ which normalizes
$\tilde K$. }

\medskip\noindent
{\it Proof of Claim:} Since the center of $G$ is connected, the map
$G(L)\to G_{ad}(L)$ is surjective. Hence we may replace $G$ by $G_{ad}$,
in which case $B(G)_{basic} = H^1(F, G)$. Hence any representative of $b$
defines an inner form of $G$ which splits over $L$.
Assume that there exists a representative $b'$ of $[b]$ in $N(L)$ as in
the claim. Then in the corresponding inner form of $G$ there exists a
maximal torus which splits over $F^{un}$ and an $F$-rational Iwahori
subgroup fixing an alcove in the apartment for this torus. Conversely, if
the inner form of $G$ corresponding to a representative of $[b]$ has this
property, then this representative normalizes this maximal torus and the
Iwahori subgroup. Now, since any inner form of $G$ contains a maximal torus which splits
over $F^{un}$ and an $F$-rational Iwahori subgroup fixing an alcove in the
apartment for this torus, such a representative must exist, which proves
the claim.

\smallskip
Let $g\in G(L)$ be such that $\tilde K=g\tilde K_0 g^{-1}$. Then
$gb'\sigma(g)^{-1}\in \tilde K_0w\tilde K_0$, where $w\in \tilde W$
normalizes $\tilde K_0$. It follows that the component of $w$ in $W_a$ is trivial, hence by Lemma
5.4,
$w\leq t_{\mu}$. If $b'=h \tilde b\sigma(h)^{-1}$ then $gh\tilde
K_0\in X_w(b\sigma)\subset X(\mu, \tilde b)_{K_0}$.\endbeweis

\medskip\noindent
The following statement yields an inequality which goes in a sense in the
opposite direction to that defining $B(G,\mu)$.

\medskip\noindent
{\bf Proposition 5.7.} {\it Let $G$ be split over $F$. Let $S$ be a
maximal split torus over $F$. Let $b\in S(L)$. Let $\tilde K_0$ denote the
Iwahori subgroup fixing an alcove in the apartment of ${\Cal B}$
corresponding to $S$. Let $w\in \tilde W$ such
$X_w(b\sigma)\neq\emptyset$, i.e., $$\exists\ g\in G(L):
g^{-1}b\sigma(g)\in \tilde K_0 w\tilde K_0\ \ .$$  Then} $$t_{\nu_b}\leq
w\ \ .$$ \noindent {\bf Proof:} Let ${\Cal A}$ be the apartment in ${\Cal
B}$ corresponding to $\tilde S=S\otimes_FL$, and let ${\bold
a}_0\subset{\Cal A}$ be the base alcove fixed by $\tilde K_0$, and let
${\bold a} =g{\bold a}_0$, for $g\in X_w( b\sigma)$. Let $\alpha$ be the
automorphism of ${\Cal B}$ induced by $b\sigma$. Then we have for the
component $w_a$ of $w$ in the affine Weyl group $$w_a={\roman{inv}}({\bold
a}, \alpha({\bold a}))\ \ .\leqno(5.5)$$ Let ${\Cal C}$ be any quartier
corresponding to the positive vector Weyl chamber, after a choice of a
special vertex of ${\bold a}_0$. Only the germ of ${\Cal C}$ will be
relevant to us, i.e.\ ${\Cal C}$ up to translation by an element of
$X_*(S_{sc})_{{\bold R}}$. Let $\varrho_{{\Cal A}, {\Cal C}}$ be the
corresponding retraction, i.e.\ $\varrho_{{\Cal A}, {\Cal C}}=
\varrho_{{\Cal A}, {\bold a}'}$ for some alcove ${\bold a}'$ far into the
quartier, comp.\ [BT1], Prop.\ 2.9.1. Then we have the following two statements.
$$\alpha\circ
\varrho_{{\Cal A}, {\Cal C}}({\bold a})= \varrho_{{\Cal A}, {\Cal C}}\circ
\alpha({\bold a})\leqno(5.6)$$

\par\noindent
(5.7)\qquad For any two alcoves ${\bold a}, {\bold a}'$ $${\roman{inv}}
(\varrho_{{\Cal A}, {\Cal C}}({\bold a}), \varrho_{{\Cal A}, {\Cal C}}
({\bold a}'))\leq {\roman{inv}}({\bold a}, {\bold a}')\ \ .$$ To see
(5.6), note that $\alpha$ preserves ${\Cal A}$ and the germ of ${\Cal C}$,
hence $$\alpha\circ \varrho_{{\Cal A}, {\Cal
C}}\circ\alpha^{-1}({\bold a})=\varrho_{\alpha({\Cal A}), \alpha({\Cal
C})}({\bold a})
=\varrho_{{\Cal A}, \alpha({\Cal C})} ({\bold a})\ \ .\leqno(5.8)$$
Since $\alpha$ also preserves the germ of ${\Cal C}$, it follows that
$\varrho_{{\Cal A}, \alpha({\Cal C})} ({\bold a})=\varrho_{{\Cal A}, {\Cal C}} ({\bold
a})$.

\par\noindent
To see (5.7), let ${\bold a}= {\bold a}_0, {\bold a}_1,\ldots, {\bold
a}_{\ell}={\bold a}'$ be a minimal gallery $\Gamma$ between ${\bold a}$
and ${\bold a}'$. This corresponds to a minimal decomposition of
$x={\roman{inv}}({\bold a}, {\bold a}')$, $$x=s_1\ldots s_{\ell}\ \
.\leqno(5.9)$$ Here $s_1,\ldots, s_{\ell}$ are the reflections around the
walls of type ${\bold a}_0\cap{\bold a}_1,\ldots, {\bold a}_{\ell -1}\cap
{\bold a}_{\ell}$ of the base simplex ([BT1], 2.3.10). The image of
$\Gamma$ under $\varrho= \varrho_{{\Cal A}, {\Cal C}}$ is a gallery
$\overline\Gamma$ between $\overline{{\bold a}} =\varrho({\bold a})$ and
$\overline{{\bold a}}'=\varrho({\bold a}')$. Furthermore, the type
$\overline s_1, \ldots, \overline s_{\ell}$ of $\overline\Gamma$ is
identical with that of $\Gamma$ ([BT1], 2.3.4.) Let us replace
$\overline\Gamma$ by a minimal gallery. Then we can write $\overline
x={\roman{inv}}(\overline{{\bold a}}, \overline{{\bold a}}')$ as
$$\overline x=s_{i_1}\ldots s_{i_k}\leqno(5.10)$$ ([BT1], 2.1.9 and
2.1.11). Hence $\overline x\leq x$, which proves (5.7).

\par\noindent
We now apply this to ${\bold a}$ and $\alpha({\bold a})$. But for any
$\overline{{\bold a}}$ contained in ${\Cal A}$ we have
$${\roman{inv}}(\overline{{\bold a}}, \alpha(\overline{{\bold a}}))=
(t_{\nu})_a\ \ .\leqno(5.11)$$ Here $\nu=\nu_b\in X_*(S)$. Hence using (5.6) and (5.7),
$$\align
(t_{\nu})_a & = {\roman{inv}}(\varrho({\bold a}), \alpha\varrho({\bold
a}))= {\roman{inv}}(\varrho({\bold a}), \varrho(\alpha({\bold a}))\\ &
\leq {\roman{inv}}({\bold a}, \alpha({\bold a}))= w_a\ \ .
\endalign$$
Taking into account the definition of the Bruhat order on $\tilde W$, the
assertion follows.\qed

\medskip\noindent
Whereas Proposition 5.6 concerned the case of a basic element $b$, the
following proposition treats the other extreme, namely, unramified
elements $b$.

\medskip\noindent
{\bf Proposition 5.8.} {\it Let $G, S, \tilde K_0$ and $b\in S(L)$ be as
in the previous proposition. Assume that $G_{der}$ is simply connected.
Then $X(\mu, b)_{K_0}\neq\emptyset\Leftrightarrow [b]\in B(G, \mu)$. }

\medskip\noindent
{\bf Proof:} If $X(\mu,  b)_{K_0}\neq \emptyset$, there exists $w\in
{\roman{Adm}}(\mu)$ such that $X_w( b\sigma) \neq\emptyset$. This implies
$\kappa(w)= \kappa(b)$ and $t_{\nu_b}\leq w$, by Lemma 5.4 and Proposition
5.7. Since $w\leq t_{\mu'}$ for some conjugate $\mu'$ of $\mu$ it follows
that $t_{\nu_b}\leq t_{\mu'}$ which implies that $W_0t_{\nu_b}W_0\leq
W_0t_{\mu}W_0$ and hence $\nu_b\mathop{\leq}\limits^! \mu$, i.e.\ $[b]\in
B(G, \mu)$.

\par
Conversely, let $[b]\in B(G, \mu)$. Hence $\nu_b$ and $\mu$ are both
dominant elements in $X_*(S)$ with $\nu_b\leq \mu$. However, for any
alcove ${\bold a}$ in the apartment ${\Cal A}$ corresponding to $S$, we
have ${\roman{inv}}({\bold a}, b\sigma({\bold a}))= (t_{\nu_b})_a$, hence
$X_{t_{\nu_b}}(b\sigma )\neq\emptyset$. Hence it suffices to see that
$t_{\nu_b}\in {\roman{Adm}}(\mu)$. But $\nu_b\leq\mu$, hence since
$G_{der}$ is simply connected, $\nu_b\mathop{\leq}\limits^!\mu$. Therefore
by Proposition 3.11 $t_{\nu_b}\leq t_{\mu}$, i.e.\ $t_{\nu_b}\in
{\roman{Adm}}(\mu)$.
\endbeweis

\medskip\noindent
As mentioned above, the sets $X(\mu,b)_K$ should have the structure of an
algebraic variety over the residue field $\F$ of ${\Cal O}_L$, at least
when $\{\mu\}$ is minuscule. For their dimension there is a conjectural
formula when $b$ is basic. To state it we first mention the following
result. The set $B(G,\mu)$ is partially ordered (a finite poset) by
$[b]\leq [b']$ iff $\overline\nu_{[b]}\leq \overline\nu_{[b']}$ in
$\Afr^+$. That this is indeed a partial order follows from the fact that
the map $(\overline\nu,\kappa): B(G)\to \Afr_+\times X^*(\hat
Z(G)^\Gamma)$ is injective [K4], [RR].

\medskip\noindent
{\bf Theorem 5.9.} (Chai [C2]) {\it Assume $G$ quasisplit over $F$.

(i) Any subset of $B(G,\mu)$ has a join, i.e.\ a supremum.

(ii) The poset $B(G,\mu)$ is ranked, i.e.\ any two maximal chains between
two comparable elements have the same length.

(iii) Let $[b], [b']\in B(G,\mu)$ with $[b]\leq [b']$. Then the length of
the maximal chain between $[b]$ and $[b']$ is given by
$${\roman{length}}([b], [b'])=\sum\limits_{i=1}^{\ell} ([ \langle
\omega_i, \overline\nu_{[b']}\rangle - \langle
\omega_i,\overline\mu^*\rangle] -[\langle \omega_i,
\overline\nu_{[b]}\rangle -\langle \omega_i, \overline\mu^*\rangle])\ .$$}
Here $\omega_1,\ldots, \omega_{\ell}$ are the fundamental $F$-weights of
the adjoint group $G_{ad}$, i.e.\ $\langle \omega_i,
\alpha_j^{\vee}\rangle =\delta_{ij}$ for any simple relative coroot
$\alpha_j^{\vee}$. Also $[x]$ denotes the greatest integer $\leq x$.

\medskip\noindent
We note that $B(G,\mu)$ has a unique minimal element, namely the unique
basic element $[b_0]$ in $B(G,\mu)$, and a unique maximal element, namely
the $\mu$-ordinary element $[b_1]=[b_{\mu}]$ for which
$\overline\nu_{[b_1]}=\overline\mu^*$. Given (i) and (ii) of Theorem 5.9,
the formula in (iii) is equivalent to $${\roman{length}}([b], [b_1])
=-\sum\limits_{i=1}^{\ell} [\langle\omega_i, \overline\nu_{[b]} \rangle
-\langle\omega_i, \overline\mu^*\rangle]\ \ . \leqno(5.12)$$ The dimension
formula for $X(\mu, b)_K$ may now be given as follows.

\medskip\noindent
{\bf Conjecture 5.10.} {\it Let $K$ be a hyperspecial maximal parahoric.
Let $[b]= [b_0]\in B(G,\mu)$ be basic. Then $X(\mu, b)_K$ is equidimensional of dimension
$$\align {\roman{dim}}\ X(\mu, b)_K &
=
\langle 2\varrho, \overline\mu^*\rangle -{\roman{length}} ([b], [b_{\mu}])
\\
&
=
\langle 2\varrho, \overline\mu^*\rangle -\sum\limits_{i=1}^{\ell} [-\langle
\omega_i, \overline\mu^*\rangle ]\ \ .
\endalign$$
Here $\varrho$ denotes the half-sum of all positive roots.}

\medskip\noindent
The motivation for this formula comes from global considerations connected
with the Newton strata of Shimura varieties, comp.\ Theorem 7.4 below. It
would be interesting to extend this conjecture to the non-basic case and
also to the case when $K$ is no longer hyperspecial.

\bigskip

\subheading{6. Relations to local models}

\bigskip

We continue with the notation of the last two sections. In particular, $G$
denotes a connected reductive group over $F$ and $\{\mu\}$ is a conjugacy
class of one-parameter subgroups of $G$. Again $E$ is the field of
definition of $\{\mu\}$. Let $K$ be a parahoric subgroup of $G(F)$ and
$\tilde K$ the corresponding parahoric subgroup of $G(L)$. We denote by
${\Cal G}={\Cal G}_K$ the group scheme over ${\Cal O}_F$ corresponding to
$K$, cf.\ Remark 2.2, (ii).

To these data one would like to associate the {\it local model,} a
projective scheme ${\Cal M}^{\roman{loc}}= {\Cal
M}^{\roman{loc}}(G,\mu)_K$ over ${\roman{Spec}}\ {\Cal O}_E$, equipped
with an action of ${\Cal G}_{{\Cal O}_E}$, at least if $\{\mu\}$ is
minuscule. It is not clear at the moment how to characterize ${\Cal
M}^{\roman{loc}}$ or how to construct it in general. It should have at
least the following properties.

\item{(i)} ${\Cal M}^{\roman{loc}}$ is flat over ${\roman{Spec}}\ {\Cal
O}_E$ with generic fiber isomorphic to $G/P_{\mu}$. Here $P_{\mu}$ denotes
the conjugacy class of parabolic subgroups corresponding to $\{\mu\}$.
\item{(ii)} There is an identification of the geometric points of the
special fiber, $${\Cal M}^{\roman{loc}}(\overline\kappa_E)= \{ g\in
G(L)/\tilde K;\ \tilde Kg\tilde K\in {\roman{Adm}}_{\tilde K}(\mu)\}\ \
.$$
\item{(iii)} ${\Cal M}^{\roman{loc}}(G,\mu)_K$ is functorial in $K$ and in
$G$.

\medskip\noindent
{\bf Examples 6.1.} (i) If $K$ is hyperspecial, then we set ${\Cal
M}^{\roman{loc}}(G,\mu)_K= {\Cal G}_{{\Cal O}_E}/{\Cal P}_{\mu}$, where
${\Cal P}_{\mu}$ is in the conjugacy class of parabolic subgroups in ${\Cal
G}_{{\Cal O}_E}$ corresponding to $\{\mu\}$. In this case ${\Cal
M}^{\roman{loc}}(G,\mu)_K$ is smooth over ${\roman{Spec}}\ {\Cal O}_E$.
Property (ii) follows from Corollary 3.12.

\medskip
\item{(ii)} Let $V$ be an $F$-vector space of dimension $n$. Let $G=GL(V)$
and let $\{\mu\}$ be minuscule of weight $r$ for some $0\leq r\leq n$,
i.e.\ $\omega_r\in\{\mu\}$, where $\omega_r(t)= {\roman{diag}}(t,\ldots,
t, 1,\ldots, 1)$ with $r$ times $t$ and $n-r$ times 1. Let $e_1,\ldots,
e_n$ be a basis of $V$ and, for $0\leq i\leq n-1$, let
$$\Lambda_i={\roman{span}}_{{\Cal O}_F}\{ \pi^{-1}e_1,\ldots, \pi^{-1}e_i,
e_{i+1},\ldots, e_n\}\ \ .$$ For a non-empty subset $I$ of $\{ 0,1,\ldots,
n-1\}$, let $K=K_I$ be the parahoric subgroup of $G(F)$ which is the
common stabilizer of the lattices $\Lambda_i$, for $i\in I$. The local
model ${\Cal M}^{\roman{loc}} = {\Cal M}^{\roman{loc}}(G,\mu)_K$ for this
triple $(G,\{\mu\}, K_I)$ represents the following moduli problem on
$(Sch/{\roman{Spec}}\ {\Cal O}_F)$ (in this case $E=F$). To $S$ the
functor associates the isomorphism classes of commutative diagrams of
${\Cal O}_S$-modules,

$$\matrix \Lambda_{i_0,S} & \longrightarrow & \Lambda_{i_1,S}
&\longrightarrow & \ldots &\longrightarrow & \Lambda_{i_m,S}
&\buildrel\pi\over\longrightarrow & \Lambda_{i_0,S}
\\
\bigcup && \bigcup &&&& \bigcup && \bigcup
\\
{\Cal F}_{i_0} &\longrightarrow & {\Cal F}_{i_1} &\longrightarrow & \ldots
&\longrightarrow & {\Cal F}_{i_m} &\longrightarrow & {\Cal F}_{i_0} & .
\endmatrix$$
Here $I=\{ i_0<i_1<\ldots <i_m\}$ and we have set
$\Lambda_{i,S}=\Lambda_i\otimes_{{\Cal O}_F}{\Cal O}_S$. It is required
that ${\Cal F}_{i_j}$ is a locally free ${\Cal O}_S$-module of rank $r$
which is locally a direct summand of $\Lambda_{i_j, S}$. The main result
of the paper [G1] of G\"ortz is that ${\Cal M}^{\roman{loc}}$ satisfies
the conditions (i) and (ii) above.

\medskip
\item{(iii)} Let $V$ be a $F$-vector space of dimension $2n$ with a
symplectic form $\langle\ ,\ \rangle$. Let $G=GSp(V,\langle\ ,\ \rangle)$
and let $\{\mu\}$ be minuscule of weight $n$. Let $e_1,\ldots, e_{2n}$ be
a symplectic basis of $V$, i.e. $$\langle e_i, e_j\rangle = 0,\ \langle
e_{i+n}, e_{j+n}\rangle =0,\ \langle e_i, e_{2n-j+1}\rangle = \delta_{ij}$$
for $i,j=1,\ldots, n$. Let $I$ be a non-empty subset of $\{ 0,\ldots,
2n-1\}$ which with $i\neq 0$ also contains $2n-i$. Let $K=K_I$ be the
parahoric subgroup of $G(F)$ which is the common stabilizer of the
lattices $\Lambda_i$ for $i\in I$. The local model ${\Cal M}^{\roman{loc}}
={\Cal M}^{\roman{loc}}(G,\mu)_K$ for the triple $(G,\{\mu\}, K_I)$
represents the moduli problem on $(Sch/{\Cal O}_F)$ which to $S$
associates the objects $({\Cal F}_{i_0},\ldots, {\Cal F}_{i_m})$ of the
local model for $(GL(V), \{\mu\}, K_I)$ as in Example (ii) above which
satisfy the following additional condition. For each $i\in I$ the
composition $${\Cal F}_i\longrightarrow \Lambda_{i, S}\simeq
\hat\Lambda_{2n-i,S}\longrightarrow \hat{\Cal F}_{2n-i}$$ is the zero map.
Here ``hat'' denotes the dual ${\Cal O}_S$-module.

By the main result of [G2], ${\Cal M}^{\roman{loc}}$ satisfies the
conditions (i) and (ii) above.

\medskip
\item{(iv)} Let $G=R_{F'/F}(GL_n)$ or $G=R_{F'/F}(GSp_{2n})$, where $F'$
is a totally ramified extension. Let $\{\mu\}$ be a minuscule conjugacy
class of one-parameter subgroups and let $K$ be a parahoric subgroup of
$G(F)$. In [PR1] resp.\ [PR2] local models ${\Cal
M}^{\roman{loc}}(G,\mu)_K$ are constructed which satisfy conditions (i)
and (ii) above. But in these cases it seems difficult to describe the
functors that these local models represent.

\medskip
In all these examples, property (ii) for local models can be considerably
strengthened by identifying the special fiber of ${\Cal M}^{\roman
loc}(G,\mu)_K$ with a closed subscheme of the partial flag variety
corresponding to $\tilde K$ of the loop group over $\kappa_E$ associated
to $G$, [G1], [PR2]. Here, to be on the safe side, we are assuming $G$
split. Via this identification there is a link between the theory of local
models and the geometric Langlands program of Beilinson, Drinfeld et al.\
[BD].

\medskip
The true significance of the local models becomes more transparent when
they appear in the global context of Shimura varieties, cf.\ (7.1). Here
we explain their relation with the sets $X(\mu, b)_K$. Let $$\tilde K_1=
{\roman{ker}}({\Cal G}({\Cal O}_L)\longrightarrow {\Cal
G}(\overline\kappa_E))\ \ .$$ Let $$X(\mu,b)_{K_1}=\{ g\in G(L) /\tilde
K_1;\ {\roman{inv}}(g\tilde K, b\sigma(g)\tilde K)\in
{\roman{Adm}}_K(\mu)\}\ \ .$$ In other words, $X(\mu, b)_{K_1}$ is the
inverse image of $X(\mu, b)_K$ under $G(L)/\tilde K_1\to G(L)/\tilde K$.
We define a map $$\tilde\gamma: X(\mu, b)_{K_1}\longrightarrow
{\Cal M}^{\roman{loc}}(G,\mu)_K(\overline\kappa_E)\leqno(6.1)$$ by
$$\tilde\gamma(g\tilde K_1)= g^{-1} b\sigma(g)\cdot \tilde K\ \ .$$ This
is well-defined since $\tilde K_1$ acts trivially on ${\Cal
M}^{\roman{loc}}(\overline\kappa_E)$. Noting that $\tilde K/\tilde K_1$ is
a principal homogeneous space under ${\Cal G}(\overline\kappa_E)$, we may
write $\tilde\gamma$ more suggestively as a map on geometric points of
algebraic stacks, $$\gamma: X(\mu, b)_K\longrightarrow [{\Cal
M}^{\roman{loc}}/{\Cal G}\otimes_{{\Cal O}_F}{\Cal
O}_E](\overline\kappa_E)\ \ .\leqno(6.2)$$ It should be possible, at least
if $\{\mu\}$ is minuscule, to equip $X(\mu, b)_K$ with the structure of an
algebraic variety over $\overline\kappa_E$ and the map $\gamma$ should be
induced by a morphism of algebraic stacks over $\overline\kappa_E$,
$$X(\mu,b)_K\longrightarrow [{\Cal M}^{\roman{loc}}\otimes_{{\Cal
O}_E}\overline\kappa_E / {\Cal G}\otimes_{{\Cal O}_F}\overline\kappa_E]\ \
.$$ Furthermore, this morphism should be compatible with Weil descent data
over $\kappa_E$ on source and target.

\par\noindent
After G\"ortz's theorems the most interesting question is the following
conjecture, comp.\ also [P], Conj.\ 2.12.

\medskip\noindent
{\bf Conjecture 6.2.} {\it Assume that $G$ is unramified over $F$. Let
${\Cal M}^{\roman{loc}}(G,\mu)_K$ be the local model over ${\roman{Spec}}\
{\Cal O}_E$, corresponding to a parahoric subgroup $K$ of $G(F)$ and a
minuscule conjugacy class of cocharacters $\{\mu\}$, with its action of
${\Cal G}\otimes_{{\Cal O}_F}{\Cal O}_E$. Then there exists a ${\Cal
G}\otimes_{{\Cal O}_F}{\Cal O}_E$-equivariant blowing up in the special
fiber $\tilde{\Cal M}^{\roman{loc}} (G,\mu)_K\to {\Cal M}^{\roman{loc}}
(G,\mu)_K$ which has semistable reduction. }

\medskip\noindent
In Example 6.1, (ii) the conjecture above is trivial for $r=1$ (in this
case ${\Cal M}^{\roman{loc}}$ has semistable reduction). For $r=2$,
Faltings [F2] has constructed an equivariant blowing-up with semistable
reduction, i.e., the conjecture holds in this case, comp.\ also [L]. In
Example 6.1, (iii) the existence of a semistable blowing-up is due to de
Jong [J] for $n=2$ and to Genestier [Ge] for $n=3$.

\bigskip

\heading {II.} {Global theory}
\endheading

\bigskip

\subheading{7. Geometry of the reduction of a Shimura variety}

\bigskip

In the global part we use the following notation. Let $(\G, \{ h\})$ be a
Shimura datum, i.e.\ $\G$ is a connected reductive group over $\Q$ and $\{
h\}$ a $G(\R)$-conjugacy class of homomorphisms from $R_{\C/\R}\G_m$ to
$\G_{\R}$ satisfying the usual axioms. We fix a prime number $p$. We let
$\K$ be an open compact subgroup of $\G(\A_f)$ which is of the form $\K=
K^p.K$, where $K^p\subset G(\A_f^p)$ and where $K=K_p$ is a parahoric
subgroup of $G(\Q_p)$. We also assume that $K^p$ is sufficiently small to
exclude torsion phenomena.

The corresponding Shimura variety $Sh(\G, \{ h\})_{\K}$ is a
quasi-projective variety defined over the Shimura field $\E$, a finite
number field contained in the field $\overline\Q$ of algebraic numbers. It
is the field of definition of the conjugacy class $\{\mu_h\}$, where
$\mu_h$ is the cocharacter corresponding to $h\in \{ h\}$. Let
$G=\G\otimes_{\Q}\Q_p$ and $F=\Q_p$. After fixing an embedding of
$\overline\Q$ into an algebraic closure of $\Q_p$, we obtain a conjugacy
class $\{\mu\}$ of one-parameter subgroups of $G$. It is defined over the
completion $E$ of $\E$ in the $p$-adic place ${\bold p}$ induced by this
embedding. We have therefore obtained by localization a triple
$(G,\{\mu\}, K)$ as in the local part relative to $F=\Q_p$.

We make the basic assumption that the Shimura variety has a good integral
model over ${\roman{Spec}}\ {\Cal O}_{\E,{\bold p}}$. Although we do not
know how to characterize it, or how to construct it in general, we know a good
number of examples all related to moduli spaces of abelian varieties. The
``facts'' stated below all refer to these moduli spaces, and the
conjectures also concern these moduli spaces or are extrapolations to the
general case. We denote by $Sh(G,h)_{\bold K}$ the model over
${\roman{Spec}}\ {\Cal O}_E$ which is obtained by base change ${\Cal
O}_{\E, {\bold p}}\to {\Cal O}_E$ from this good integral model.

The significance of the local model ${\Cal M}^{\roman{loc}}(G,\mu)_{K}$ is
given by the relatively representable morphism of algebraic stacks over
${\roman{Spec}}\ {\Cal O}_E$, $$\lambda: Sh(G,h)_{\bold K}\longrightarrow
[ {\Cal M}^{\roman{loc}} (G,\mu)_{K} / {\Cal G}_{{\Cal O}_E}]\ \
.\leqno(7.1)$$ Here ${\Cal G}$ is the group scheme over ${\roman{Spec}}\
\Z_p$ corresponding to $K_p$. The morphism $\lambda$ is of relative dimension
${\roman{dim}}\ G$. This statement is proved in those cases where an
integral model of the Shimura variety $Sh(\G, h)_{\bold K}$ exists [PR2], [RZ2].

We now turn to the special fiber. By associating to each point of
$Sh(G,h)_{\bold K}(\overline\kappa_E)$ the isomorphism class of its
rational Dieudonn\'e module (again this makes sense only for those Shimura
varieties which are moduli spaces of abelian varieties), we obtain a map
$$\delta : Sh(G,h)_{\bold K}(\overline\kappa_E)\longrightarrow B(G)\ \ .$$
We note that by Proposition 4.4 the image of $\delta$ is contained in
$B(G,\mu)$, provided that $G$ is unramified.

\medskip\noindent
{\bf Conjecture 7.1.} ${\roman{Im}}(\delta)= B(G,\mu)$.

\medskip\noindent
In particular, we expect that the basic locus is non-empty. Let
$\overline{Sh(G,h)}_{\bold K} = Sh(G,h)_{\bold K}\otimes_{{\Cal
O}_E}\kappa_E$. The {\it basic locus} $\overline{Sh(G,h)}_{{\bold K},\
{\roman{basic}}}$ is the set of points whose image under $\delta$ is the
unique basic element $[b_0]$ of $B(G,\mu)$. More generally, for $[b]\in
B(G)$, let $${\Cal S}_{[b]}=\delta^{-1}([b])\ \ .\leqno(7.2)$$

\medskip\noindent
{\bf Proposition 7.2.} ([RR]) {\it Each ${\Cal S}_{[b]}$ is a locally
closed subvariety of $\overline{Sh(G,h)}_{\bold K}$. Furthermore, for
$[b], [b']\in B(G)$, we have $${\Cal S}_{[b]}\cap{\roman{closure}}({\Cal
S}_{[b']})\neq\emptyset\Longrightarrow [b]\leq [b']\ \ .$$}

This is the group-theoretic version of Grothendieck's semicontinuity
theorem, according to which the Newton vector of an isocrystal decreases
under specialization (in the natural partial order on $(\R^n)_+$, cf.\
(4.5)). The subvarieties ${\Cal S}_{[b]}$ are called the {\it Newton
strata of} $\overline{Sh(G,h)}_{\bold K}$. The Newton stratification of
the special fiber is very mysterious.

\medskip\noindent
{\bf Questions 7.3.} Assume that $K_p$ is hyperspecial. Let $[b], [b']\in
B(G,\mu)$ with $[b]\leq [b']$.
\item{(i)} Is ${\Cal S}_{[b]}\cap{\roman{closure}}({\Cal
S}_{[b']})\neq\emptyset$?
\item{(ii)} Is ${\Cal S}_{[b]}\cap{\roman{closure}}({\Cal
C})\neq\emptyset$, for every irreducible component ${\Cal C}$ of ${\Cal
S}_{[b']}$?
\item{(iii)} Is ${\Cal S}_{[b]}\subset {\roman{closure}}({\Cal
S}_{[b']})$?

\medskip\noindent
Obviously, (i) is implied by (iii). Here part (iii) has become known as
the strong Grothendieck conjecture and (i) as the weak Grothendieck
conjecture, although this denomination is somewhat abusive.

\medskip\noindent
{\bf Theorem 7.4.} (Oort [O]) {\it For the Shimura variety associated to
$GSp_{2n}$ {\it (the Siegel moduli space),} question 7.3 (iii) has an
affirmative answer. Also in this case, each Newton stratum ${\Cal
S}_{[b]}$ is equidimensional of codimension in $\overline{Sh(G,h)}_{\bold
K}$ equal to $${\roman{codim}}\ {\Cal S}_{[b]}= {\roman{length}}([b],
[b_{\mu}])\ \ ,$$ comp.\ Theorem 5.9. Here $[b_\mu]=[b_1]$ denotes the
$\mu$-ordinary element of $B(G,\mu)$, cf.\ (5.12).}

\medskip\noindent
{\bf Conjecture 7.5.} (Chai [C2]) {\it Assume $K_p$ hyperspecial. Each
Newton stratum ${\Cal S}_{[b]}$ is equidimensional of codimension given by
the above formula.}

Since the basic element in $B(G,\mu)$ is minimal, the basic locus
$\overline{Sh(G,h)}_{K,\ {\roman{basic}}}$ is a closed subvariety of the
special fiber. This variety has been studied in many cases ([LO], [Ka], [Ri]); it
is conceivable that one can give a group-theoretical ``synthetic''
description of it in general. At the other extreme is the $\mu$-ordinary
element $[b_{\mu}]\in B(G,\mu)$. It is the unique maximal element of
$B(G,\mu)$, cf.\ (5.12).

\medskip\noindent
{\bf Conjecture 7.6.} (Chai) {\it Let $K_p$ be hyperspecial. The orbit of
any point of ${\Cal S}_{[b_{\mu}]}$ under ${\bold G}(\A_f^p)$ is dense in
$\overline{Sh(G,h)}_{\bold K}$.}

Here the action of ${\bold G}(\A_f^p)$ is via Hecke correspondences. In
this direction we have the following results.

\medskip\noindent
{\bf Theorem 7.7.} (Chai [C1]) {\it The conjecture 7.6 is true for the
Siegel moduli space.}

\medskip\noindent
{\bf Theorem 7.8.} (Wedhorn [W]) {\it We assume that the Shimura variety
$Sh(G,h)_{\bold K}$ corresponds to a PEL-moduli problem of abelian
varieties. Let $K_p$ be hyperspecial. The $\mu$-ordinary stratum ${\Cal
S}_{[b_{\mu}]}$ is dense in $\overline{Sh(G,h)}_{\bold K}$. }

\medskip\noindent
The hypothesis that $K_p$ be hyperspecial in Wedhorn's theorem is indeed
necessary, as the examples of Stamm [S] relative to the
Hilbert-Blumenthal surfaces with Iwahori level structure at $p$ and of
Drinfeld [D] relative to a group which is ramified at $p$ show.

We finally relate the maps $\gamma$ and $\lambda$. Let $[b]\in B(G,\mu)$
and let $b\in G(L)$ be a representative of $[b]$. As in Definition 4.1, we
let $$J_b(\Q_p)= \{ g\in G(L);\ g^{-1} b \sigma(g)= b\}\ \ .$$
The Newton stratum ${\Cal S}_{[b]}$ has a covering $\tilde{\Cal S}_{b}$,
for which we fix an isomorphism of the isocrystal in the variable point
$x\in {\Cal S}_{[b]}$ with the model isocrystal with $G$-structure
determined by $b$. Then $\tilde{\Cal S}_b$ is a principal homogeneous
space under $J_b(\Q_p)$ over ${\Cal S}_{[b]}$. The relation between
$\gamma$ and $\lambda$ is then given by a commutative diagram of morphisms
of algebraic stacks (over $\overline\kappa_E$, or even compatible with
Weil descent data over $\kappa_E$), $$\matrix \tilde{\Cal S}_b & \buildrel
\tilde\lambda\over\longrightarrow & X(\mu, b)_{K_p}
\\
\big\downarrow && \big\downarrow\rlap{$\scriptstyle\gamma$}
\\
{\Cal S}_{[b]} & \buildrel\lambda\over \longrightarrow & \left[
\overline{\Cal M}^{\roman{loc}}(G,\mu)_{K_p} /\overline{\Cal
G}_{\kappa_E}\right] & .
\endmatrix$$
Here $\overline{\Cal M}^{\roman{loc}}$ resp.\ $\overline{\Cal G}$ denote
the special fibers of the local model resp.\ the group scheme
corresponding to $K_p$, and by $\lambda$ we denoted the restriction of
(7.1) to ${\Cal S}_{[b]}$.

\bigskip

\subheading{8. Pseudomotivic and quasi-pseudomotivic Galois gerbs}

\bigskip

In this section and the next one we wish to give a conjectural description
of the point set of $Sh(G,h)_{\bold K} (\bar\kappa_E)$ with its action of
the Frobenius automorphism. This description is modeled on the one given
in [LR], but differs from it in an important detail, comp.\ Remark 9.3.
The idea is to partition the point set into ``isogeny classes'', as was
done in the case of the elliptic modular curve in (1.1), and then to
describe the point set of the individual isogeny classes in a manner
reminiscent of (1.1) and (1.2) resp.\ (1.9) in the elliptic modular case.
According to an idea of Grothendieck, the set of isogeny classes will be
described in terms of representations of certain Galois gerbs. In this
section we introduce these Galois gerbs, following Reimann's book [Re1].

We first explain our terminology concerning Galois gerbs.
 Let $k$ be a field of characteristic zero and let $\Gamma = Gal (\bar k/k)$ be
  the Galois group of a chosen algebraic closure. A {\it Galois gerb over} $k$ is an
  extension of topological groups
$$ 1 \to G (\bar k) \to {\Cal G} \buildrel q \over \longrightarrow \Gamma
\to 1.\leqno(8.1) $$
Here $G$ denotes a linear algebraic group over $\bar k$ and is called the
{\it kernel} of ${\Cal G}$.
 The topology on $\Gamma$ is the Krull topology and the topology on $G(\bar k)$ is the
 discrete topology. The extension ${\Cal G}$ is required to satisfy the following two conditions:

(i) For any representative $g_\sigma \in {\Cal G}$ of $\sigma \in \Gamma$,
the automorphism $g \mapsto g_\sigma g g_\sigma^{-1}$ of $G(\bar k)$ is a
$\sigma$-linear algebraic
 automorphism.

(ii) Let $K/k$ be a finite extension over which $G$ is defined. Let
$\Gamma_K = Gal (\bar k/K)$ be the corresponding subgroup of $\Gamma$. We
choose a section of ${\Cal G} \to \Gamma$ over $\Gamma_K$ such that the
automorphism
$$ g \longmapsto g_\sigma g g_\sigma^{-1}, \quad g \in G(\bar k) $$
defines the $K$-structure on $G$. Then the resulting bijection
$$ q^{-1} (\Gamma_K) = G(\bar k) \times \Gamma_K $$
is a homomorphism.

A {\it morphism between Galois gerbs} $\varphi : {\Cal G} \to {\Cal G}'$
is a continuous map of extensions which induces the identity map on
$\Gamma$ and an algebraic homomorphism on the kernel groups. Two
homomorphisms $\varphi_1$ and $\varphi_2$ are called {\it equivalent} if
there exists $g' \in G' (\bar k)$ with $\varphi_2 = {\roman{Int}} (g)
\circ \varphi_1$. A {\it neutral gerb} is one isomorphic to the
semi-direct product
$$ {\Cal G}_G = G (\bar k) \rtimes \Gamma $$
associated to an algebraic group $G$ over $k$. Sometimes we will have to
consider a slightly more general notion. Let $k'$ be a Galois extension of
$k$ contained in $\bar k$ with Galois group $\Gamma' = Gal (k'/k)$. Then
one defines in the obvious way the notion of a $k'/k$-Galois gerb, which
is an extension
$$ 1 \to G(k') \to {\Cal G} \to \Gamma' \to 1,\leqno(8.2) $$
where $G$ is an algebraic group defined over $k'$. A $k'/k$-Galois gerb
defines in the obvious way a $\bar k/k$-Galois gerb (i.e. a Galois gerb):
one first pulls back the extension by the surjection $\Gamma \to \Gamma'$
and then pushes out via $G(k') \to G(\bar k) = (G \otimes_{k'} \bar k)
(\bar k)$.

In the sequel we will have to deal with projective limits of Galois gerbs
of the previous kind. We transpose the above terminology to them. In
particular two morphisms of pro-Galois gerbs will be called equivalent if
they are projective limits of equivalent morphisms of Galois gerbs (in
[Re1], B.1.1, this is called {\it algebraically equivalent}).

An important example is given by the Dieudonn\'e gerb over ${\Q}_p$, cf.
[Re1], B.1.2. For every $n \in {\Z}, \; n \ge 1$, there is an explicitly
defined ${\Q}_p^{un}/{\Q}_p$-gerb ${\Cal D}_n$ with kernel group ${\G}_m$.
For $n'$ divisible by $n$ there is a natural homomorphism ${\Cal D}_{n'}
\to {\Cal D}_n$ inducing the map $x \mapsto x^{n'/n}$ on the kernel
groups. Let
$${\Cal D}_0=\lim\limits_{\longleftarrow} {\Cal D}_n\leqno(8.3)$$
be the pro-${\Q}_p^{un}/{\Q}_p$-Galois gerb defined by this projective
system. Then in ${\Cal D}_0$ there is an explicit representative
$d_\sigma$ of the Frobenius element. The {\it Dieudonn\'e gerb} ${\Cal D}$ is
the pro-$\bar {\Q}_p/{\Q}_p$-Galois gerb defined by ${\Cal D}_0$.

Another Galois gerb of relevance to us is the {\it weight gerb} ${\Cal
W}$. This is the Galois gerb over $\R$ with kernel $\G_m$ which is defined
by the fundamental cocycle of ${\roman{Gal}}(\C/\R)$
($w_{\varrho,\sigma}=-1$ if $\varrho=\sigma=$ complex conjugation;
otherwise $w_{\varrho,\sigma}=1$).

We now recall some pertinent facts about the pro-Galois gerbs appearing in
the title of this section. We fix an algebraic closure $\bar {\Q}$ of
${\Q}$ and for every place $\ell$ of ${\Q}$ an embedding $\bar{\Q} \subset
\bar {\Q}_\ell$. Let $L/{\Q}$ be a finite Galois extension contained in
$\bar {\Q}$.

There is an initial object $(Q^L, \nu (\infty)^L, \nu(p)^L)$ in the
category of all triples $(T, \nu_\infty, \nu_p)$ where $T$ is a
${\Q}$-torus which splits over $L$ and such that $\nu_\infty, \nu_p \in
X_\ast (T)$ satisfy
$$ [L_\infty : {\R}]^{-1} \cdot Tr_{L/{\Q}} (\nu_\infty) + [L_p :
{\Q}_p]^{-1} \cdot Tr_{L/{\Q}} (\nu_p)=0,\leqno(8.4) $$
cf. [Re1], B.2.2.

Similarly, assume that $L$ is a $CM$-field and denote by $L_0$ its maximal
totally real subfield. Then there is an initial object $(P^L, \nu(\infty)^L,
\nu(p)^L)$ in the category of all triples $(T, \nu_\infty, \nu_p)$ where
$T$ is a ${\Q}$-torus which splits over $L$ and such that $\nu_\infty,
\nu_p \in X_\ast (T)$ are defined over ${\Q}$ and ${\Q}_p$ respectively
and such that
$$ \nu_\infty + [L_p : {\Q}_p]^{-1} \cdot Tr_{L/L_0} (\nu_p) =
0,\leqno(8.5) $$
cf. [Re1], B.2.3. Since obviously condition (8.5) implies condition (8.4),
there is a canonical morphism
$$(Q^L, \nu(\infty)^L, \nu(p)^L)\longrightarrow (P^L, \nu(\infty)^L,
\nu(p)^L)\ \ .\leqno(8.6)$$
If $L \subset L'$ then we obtain morphisms of tori in the opposite
direction,
$$ Q^{L'} \longrightarrow Q^L, \quad P^{L'} \longrightarrow
P^L.\leqno(8.7) $$
Let $Q$ resp. $P$ denote the pro-torus defined by this projective system.
Then there are homomorphisms of pro-tori over ${\Q}$,
$$ \nu (\infty) : {\G}_{m,{\R}} \longrightarrow Q_{\R}, \ \ \hbox{resp.}\
\nu (p): {\Dbf} \longrightarrow Q_{{\Q}_p}\leqno(8.8) $$
whose composite with $Q \to Q^L$ is $\nu (\infty)^L$ if $L$ is totally
imaginary resp. is $[L_p : {\Q}_p] \cdot \nu (p)^L$. Here $\Dbf$ denotes
the pro-torus with character group equal to $\Q$.

Similarly, we obtain
$$ \nu (\infty) : {\G}_{m} \longrightarrow P_{\bold R}, \quad \nu (p):
{\Dbf} \longrightarrow P_{{\Q}_p}.\leqno(8.9) $$
We can now introduce the pro-Galois gerbs which will be relevant for the
theory of Shimura varieties.

A {\it quasi-pseudomotivic Galois gerb} is a pro-Galois gerb $\frak{Q}$
over ${\Q}$ with kernel $Q$ together with morphisms
$$\align \zeta_\infty : {\Cal W} &\longrightarrow \frak{Q}_{\R}\\ \zeta_p
: {\Cal D} &\longrightarrow \frak{Q}_{{\Q}_p}\tag8.10\\ \zeta_\ell :
\Gamma_\ell &\longrightarrow \frak{Q}_{{\Q}_\ell}, \quad \ell \not=
\infty, p
\endalign
$$
such that $\nu(\infty)$ is induced by $\zeta_\infty$ resp. $\nu(p)$ is
induced by $\zeta_p$ on the kernel group ${\G}_{m_{\R}}$ of the weight
gerb ${\Cal W}$ resp. the kernel group ${\Dbf}$ of the Dieudonn\'e gerb
${\Cal D}$. In addition, a coherence condition on the family $\{
\zeta_\ell; \; \ell \not= \infty, p\}$ is imposed, cf. [Re1], B.2.7.
Similarly one defines a {\it pseudomotivic Galois gerb} $({\Pfr},
\zeta^P_\ell$). These pro-Galois gerbs are uniquely defined up to an
isomorphism preserving the morphisms $\zeta_\ell$ up to equivalence for
$\ell=\infty, p$ and $\ell\neq p$. Furthermore, these isomorphisms are
unique up to equivalence. There is a morphism
$$ \Qfr \longrightarrow \Pfr\leqno(8.11) $$
compatible with the morphisms $\zeta_\ell$ resp. $\zeta_\ell^P$ and
inducing the homomorphism (8.6) above on the kernel groups.

For a pair $(T,\mu)$ consisting of a ${\Q}$-torus $T$ and an element $\mu
\in X_\ast (T)$, there is associated a morphism of Galois gerbs
$$ \psi_\mu : \frak{Q} \longrightarrow {\Cal G}_T,\leqno(8.12) $$
cf. [Re1], B.2.10. This morphism factors through $\Pfr$ if and only if the
following two conditions are satisfied,

\item{(i)} the image of $\nu(\infty)$ in $X_\ast (T)$ is defined over ${\Q}$
\item{(ii)} the image of $\nu(p)$ in $X_\ast (T)\otimes {\Q}$ satisfies the Serre
condition,i.e., it is defined over a CM-field and its weight is defined
over $\Q$,
cf. [Re1], B.2.11. This is the case if $T$ itself satisfies the Serre
condition, i.e., $({\roman{id}}+\iota)(\tau-{\roman{id}})=
(\tau-{\roman{id}})({\roman{id}}+\iota)$ in ${\roman{End}}(X_*(T))$, where
$\iota$ denotes the complex conjugation and
$\tau\in{\roman{Gal}}(\overline\Q /\Q)$ is arbitrary.

\medskip\noindent
{\bf Remark 8.1.} The pseudomotivic Galois gerb was introduced in [LR]
with the aim of describing the points in the reduction of a Shimura
variety when $(G,\{ h\})$ satisfies the Serre condition. When this last
condition is dropped, the pseudomotivic Galois gerb cannot suffice for
this purpose. However, the quasi-pseudomotivic Galois gerb in [LR], introduced there to
cover the cases when the Serre condition fails, does
not exist (there is a fatal error in the construction of loc.\ cit.). Two
replacements have been suggested, one by Pfau [Pf] and one by Reimann [Re1].
We follow here the latter.

\bigskip

\subheading{9. Description of the point set in the reduction}

\bigskip

In this section we return to the notation used in section 7.  Therefore
$Sh(G, h)_{\bold K}$ is an integral model over ${\roman{Spec}}\ {\Cal
O}_E$ of the Shimura variety associated to $({\bold G}, \{h\}, {\bold
K}=K^p.K_p)$, and $\{\mu_h\}$ is the associated conjugacy class of
cocharacters of the reductive group ${\bold G}$ over ${\roman{Spec}}\
{\bold Q}$. Our purpose is to describe the set $Sh(G,h)_{\bold
K}(\overline\kappa_E)$ of our model over ${\roman{Spec}}\ {\Cal O}_E$ of
the Shimura variety $Sh(\G, h)_{\bold K}$. We make the blanket assumption
that the derived group of $\G$ is simply connected.

The description of the points in the reduction will be in terms of admissible
 morphisms of pro-Galois gerbs $$ \varphi : {\frak {Q}}
\longrightarrow {\Cal G}_G.\leqno(9.1) $$

\medskip\noindent
{\bf Definition 9.1.} A morphism (9.1) of pro-Galois gerbs over ${\bold
Q}$ is called {\it admissible} if it satisfies the four conditions a)--d)
below.

Let $D = G/G_{\roman{der}}$ and let $\mu_D$ be the image of $\{\mu_h\}$ in
$X_\ast (D)$. The first condition is global:

\medskip
a) {\it The composition ${\frak {Q}} \to {\frak {G}}_{\G} \to {\frak
{G}}_D$ is equivalent to $\psi_{\mu_D}$, cf. (8.11).}

\medskip
The next three conditions will be local, one for each place of ${\Q}$. To
formulate the next condition we remark that for $h \in \{h\}$ with
corresponding weight
 homomorphism $w_h:\G_{m,\R}\to\G_{\R}$ the pair $(w_h , \mu_h (-1))$
 corresponds to a morphism of Galois gerbs
 over ${\R}$,
$$ \xi_\infty : {\Cal W} \longrightarrow
{\Cal G}_{\G_{\R}}.\leqno(9.2) $$ b) {\it The composition
$\varphi\circ\zeta_{\infty}$ is equivalent to $\xi_\infty$.

\smallskip
{\rm c)} For any $\ell \not= \infty, p$ the composition $\varphi\circ
\zeta_{\ell}$ is equivalent to the canonical section $\xi_\ell$ of
${\Cal G}_{{\G}_{{\Q}_\ell}}$.}

\medskip
For the final condition we remark that (the equivalence class of) the
composition $\varphi\circ \zeta_p  : {\Cal D} \to {\Cal G}_{{\G_\Q}_p}$
defines an element $[b] = [b (\varphi_p)] = [b(\varphi)]$ of $B(G)$. More
precisely, let ${\Cal D}_0$ be the explicit unramified version of the
Dieudonn\'e gerb as in [Re1], B.2, comp. (8.3). Then there exists a
morphism $$ \theta_0 : {\Cal D}_0 \longrightarrow G ({\Q}_p^{un}) \rtimes
\hat{\Z}\leqno(9.3) $$ such that $\varphi\circ\zeta_p$ is equivalent to
the pullback $\bar\theta_0$ of $\theta_0$ to $\Gamma$. Then $[b]$ is the
class of $b = \theta_0 (d_\sigma)$, where $d_\sigma \in {\Cal D}_0$ is the
explicit representative of the Frobenius $\sigma$.

\medskip
d) {\it The element $[b]$ lies in }$B (G, \mu)$.

\medskip

We note that, whereas the local components $\varphi\circ \zeta_\infty$ and
$\varphi\circ \zeta_{\ell}$ $(\ell\neq p)$ are uniquely determined up to
equivalence by the Shimura data, the $p$-component $\varphi\circ \zeta_p$
is allowed to vary over a finite set of equivalence classes.

To every admissible morphism $\varphi$ we shall associate a set
$S(\varphi)$ with an action from the right of $Z({\Q}_p) \times
G({\A}^p_f)$ and a commuting action of an automorphism $\Phi$. Here $Z$
denotes the center of $G$. For $\ell \not= \infty,p$ let $$ X_\ell = \{ g
\in G(\bar{\Q}_\ell); \quad Int(g) \circ \xi_\ell =\varphi\circ
\zeta_{\ell}\}.\leqno(9.4) $$ By condition c) this set is non-empty. We
put $$ X^p = \prod\nolimits_{\ell \not= \infty,\;p}' \; X_\ell,\leqno(9.5)
$$ where the restricted product is explained in [Re1], B.3.6. The group
${\bold G}({\A}^p_f)$ acts simply transitively on $X^p$. We also put $X_p
= X (\mu, b)_{K_p}$ in the notation of (5.1), where $b\in G(L)$ is as
above. It is equipped with commuting actions of $Z({\Q}_p)$ and an
operator $\Phi$ (cf. (5.2)). Finally we introduce the group of
automorphisms $I_\varphi = Aut\,(\varphi)$. The group $I_\varphi ({\Q})$
obviously operates on $X^p$. Let $g_p \in G(\bar{\Q}_p)$ be such that $$
\varphi_p \circ \zeta_p = Int\, g_p \circ \bar\theta_0\leqno(9.6) $$ where
the notation is as in the formulation of condition d) above. Then we
obtain an embedding $$ I_\varphi ({\Q}) \longrightarrow J_b ({\Q}_p),
\quad h \mapsto g_p \,h \,g_p^{-1}\leqno(9.7) $$ where $J_b({\Q}_p)$ is
the group associated to the element $b$ which acts on $X_p = X
(\mu,b)_{K_p}$, cf. Definition 4.1. We now define $$ S(\varphi)_{K_p} =
{\lim_{\longleftarrow\atop K^p}} \ I_\varphi ({\Q}) \backslash X_p \times
X^p/K^p,\quad \hbox{ resp.}\quad S(\varphi)_{\bold K} = I_\varphi ({\Q})
 \backslash X_p \times X^p /K^p,\leqno(9.8)
$$ where the limit is over all open compact subgroups $K^p \subset {\bold
G}({\A}^p_f)$.
 On $S(\varphi)_{K_p}$ we have commuting actions of the automorphism $\Phi$
 and of $Z ({\Q}_p) \times {\bold G}({\A}^p_f)$ from the right.

\medskip\noindent
{\bf Conjecture 9.2.} {\it Assume that the derived group of ${\bold G}$ is
simply connected. Then for every sufficiently small $K^p$ there is a model
$Sh (G, \{h\})_{\bold K}$ of $Sh (\G, \{h\})_{\bold K}$ over
 ${\roman{Spec}}\, {\Cal O}_{{\bold E}_{({\bold p})}}$ such that the point set of its special fiber
  is a disjoint sum of subsets
 invariant under the action of the Frobenius automorphism over $\kappa_E$ and
 of $Z({\Q}_p)$ and $\G({\A}^p_f)$,
$$ Sh (G, h)_{\bold K} (\overline\kappa_E) = \coprod_\varphi Sh
 (G, h)_{{\bold K}, \varphi},
$$ and for each $\varphi$ a bijection $$ Sh (G, h)_{{\bold K}, \varphi} =
S(\varphi)_{\bold K}, $$ which carries the action of the Frobenius
automorphism over $\kappa_E$ on the left into the action of $\Phi$ on the
right and which commutes with the actions of $Z({\Q}_p)$ and
 of $\G({\A}^p_f)$ (for variable $K^p$) on both sides. Here the disjoint union is
 taken over a set of representatives of equivalence classes of admissible morphisms
 $\varphi : {\frak{Q}} \to {\Cal G}_{\G}$.
}

We remark that if $D$ splits over a $CM$-field and the weight homomorphism
$w_h$ is defined over ${\Q}$, every admissible morphism $\varphi :
{\frak{Q}} \to {\Cal G}_{\G}$ factors through ${\frak{P}}$, cf. [Re1],
B.3.9.

Note that we are not proposing a characterization of the model
$Sh(G,h)_{\bold K}$. In the case where $K_p$ is hyperspecial, such a
characterization was suggested by Milne [M2]. In this case we expect $Sh(G,h)_{\bold
K}$ to be smooth  over ${\roman{Spec}}\, {\Cal O}_{{\bold E}_{({\bold
p})}}$. In [Re3], Reimann gives a wider class of parahoric subgroups $K_p$
for which one should expect the smoothness of this model, and he
conjectures that this class is exhaustive.

Conjecture 9.2 has been proved by Reimann [Re1], Prop.\ 6.10 (and
Remark 4.9) and Prop.\ 7.7, in the case when $K_p$ is a maximal compact
subgroup of $G({\bold Q}_p)$. Here $\G$ is the multiplicative group of a
quaternion algebra over a totally real field in which $p$ is unramified,
which is either totally indefinite or which is unramified at all primes
above $p$. It has been proved for Shimura varieties of PEL-type by Milne
[M1], when $K_p$ is hyperspecial. In [LR] it is shown how the conjecture is
related to a hypothetical good theory of motives, comp.\ also [M3].

\medskip
\noindent {\bf Remark 9.3.} We note that if Conjecture 5.2 holds, then
each summand in Conjecture 9.2 is non-empty. In [LR] (apart from a very
special case of bad reduction) it was assumed that $K_p$ is hyperspecial,
and the admissibility condition d) was replaced by the condition that
$X(\mu, b)_{K_p}$ be non-empty. From Remark 5.3 it follows that then
$[b]\in B(G,\mu)$, i.e.\ condition d) above holds.

\medskip\noindent
{\bf Remark 9.4.} Assume Conjecture 9.2. In [RZ2] it was shown that in
certain very rare cases the Shimura variety $Sh(G, \{ h\})_{\bold K}$
admits a $p$-adic uniformization by (products of) Drinfeld upper half
spaces. The proof in loc.cit.\ is a generalization of Drinfeld's proof [D]
of Cherednik's uniformization theorem in dimension one. From the proof in
[RZ2] it is clear that this can occur only when all admissible morphisms
are locally equivalent, provided that all summands in Conjecture 9.2 are
non-empty. It comes to the same to ask that $\varphi_p\circ \zeta_p$ is
basic for any admissible morphism $\varphi$. In [K4] it is shown that
when $G$ is adjoint simple such that $B(G,\mu)$ consists of a single
element (which is then basic), then $(G,\mu)$ is the adjoint pair
associated to $(D_{1\over n}^\times, (1,0,\ldots,0))$ or $(D_{-{1\over
n}}^\times, (1,\ldots, 1,0))$, where $D_{1\over n}^\times$ resp.
$D_{-{1\over n}}^\times$ denotes the inner form of $GL_n$ associated to
the central division algebra of invariant $1/n$ resp.\ $-1/n$. In other
words, this result of Kottwitz implies in conjunction with Conjecture 9.2
and Conjecture 5.2 that there is no hope of finding cases of $p$-adic
uniformization essentially different from those in [RZ2]. In particular, in
all cases of $p$-adic uniformization the uniformizing space will be a
product of Drinfeld upper half spaces.

\bigskip

\subheading{10. The semi-simple zeta function}

\bigskip

We continue with the notation of the previous section. One ultimate goal
of the considerations of the previous section is to determine the local
factor of the zeta function of $Sh(\G,h)_{\bold K}$ at $p$. Our present
approach is through the determination of the local {\it semi-simple} zeta
function [R2]. We refer to [HN2], \S 3.1 for a systematic exposition of the
concepts of the semi-simple zeta function and semi-simple trace of
Frobenius. The decisive property of the semi-simple trace of Frobenius on
representations of the local Galois group is that it factors through the
Grothendieck group. In the case of good reduction the semi-simple zeta
function coincides with the usual zeta function.

To calculate the semi-simple zeta function we may use the Lefschetz fixed
point formula. Let $\kappa_E^n$ be the extension of degree $n$ of
$\kappa_E$ contained in $\overline\kappa_E$. For $x\in
\overline{Sh(G,h)}_K(\kappa_E^n)$ we introduce the semisimple trace
$${\roman{Contr}}_n(x)=
{\roman{tr^{ss}}}(Fr_n;\ R\Psi_x(\overline\Q_{\ell}))\ \ .\leqno(10.1)$$
Here $R\Psi(\overline\Q_\ell)$ denotes the complex of nearby cycles. By
$Fr_n$ we denote the geometric Frobenius in
$\Gal(\overline\kappa_E/\kappa_E^n)$. This is the contribution of $x$ to
the Lefschetz fixed point formula over $\kappa_E^n$. In the case of good
reduction, or more generally if $x$ is a smooth point of $Sh(G,h)_{\bold K}$, then
${\roman{Contr}}_n(x)=1$.

\par
For an admissible homomorphism $\varphi:\Qfr\to{\Cal G}_{\G}$ as in
Conjecture 9.2, we introduce the {\it contribution of} $\varphi$ (or its
equivalence class) {\it to the Lefschetz fixed point formula over}
$\kappa_E^n$, $${\roman{Contr}}_n(\varphi)=\sum_{x\in Sh(G,h)_{{\bold
K},\varphi}(\kappa_E^n)} {\roman{Contr}}_n(x)\ \ .\leqno(10.2)$$

\medskip\noindent
{\bf Definition 10.1.} {\it A morphism $\varphi:\Qfr\to{\Cal G}_{\G}$ is called
{\rm special} if there exists a maximal torus $T\subset G$ and an element
$\mu\in X_*(T)$ which defines a one-parameter subgroup of $G$ in the
conjugacy class $\{\mu_h\}$ such that $\varphi$ is equivalent to $i\circ
\psi_{\mu}$, cf.\ (8.12). Here $i:{\Cal G}_T\to {\Cal G}_{\G}$ denotes the canonical
morphism defined by the inclusion of $T$ in $G$.}

\medskip\noindent
If $K_p$ is hyperspecial (and $G_{\roman{der}}$ is simply connected, as is
assumed throughout this section), then every admissible morphism is
special ([LR], Thm.\ 5.3), at least if it factors through $\Pfr$.

\medskip\noindent
{\bf Conjecture 10.2.} {\it We have ${\roman{Contr}}_n(\varphi)=0$ unless
$\varphi$ is special.}

\medskip\noindent
For some cases of this conjecture related to $GL_2$, comp.\ [R1] and
[Re1]. Note that this is really a conjecture about bad reduction. In the
case of good reduction the cancellation phenomenon predicted by Conjecture
10.2 cannot occur since each point $x$ in $Sh(G,h)_{\bold K}(\kappa_E^n)$
contributes 1 in this case. This is compatible with the remark immediately
preceding the statement of the conjecture, which says that the conjecture
is empty if $K_p$ is hyperspecial.

\medskip
Let us explain how one would like to give a group-theoretic expression for
${\roman{Contr}}_n(x)$. Let $x$ be represented by $(x_p, x^p)\in X_p\times
X^p/K^p$ under the bijection (9.8). Let $n'=n\cdot r$, where as shortly
after (5.1) $r=[\kappa_E:\F_p]$. Since $x$ is fixed under the $n$-th power
of the Frobenius over $\kappa_E$, it is fixed under the $n'$-th power of
the absolute Frobenius and we obtain an equation of the form
$$(\Phi^{n'}x_p, x^p)=h\cdot (x_p, x^p)\ \ ,\leqno(10.3)$$ for some $h\in
I_\varphi(\Q)$. By [K1], Lemma 1.4.9, it follows that there exists $c\in
G(L)$ such that $$c\cdot h^{-1}\cdot \Phi^{n'}\cdot c^{-1}=\sigma^{n'}\ \
.\leqno(10.4)$$ This is an identity in the semi-direct product
$G(L)\rtimes \langle\sigma\rangle$, where $L$ is the completion of the
maximal unramified extension $\Q_p^{un}$ of $\Q_p$. Let $\Q_{p^{n'}}$ be
the fixed field of $\sigma^{n'}$. Then the element $\delta\in G(L)$ which
is defined by the equation $$c\cdot (b\sigma)\cdot
c^{-1}=\delta\sigma\leqno(10.5)$$ lies in $G(\Q_{p^{n'}})$. Here $b\in
G(L)$ is the element defined before (9.3). Also, always by [K1], we have
that $x'_p=c\cdot x_p$ lies in $G(\Q_{p^{n'}}) / \tilde K^{\langle
\sigma^{n'}\rangle}$. To simplify the notation put $K_{p^{n'}} = \tilde
K_p^{\langle \sigma^{n'}\rangle}$. Let $${\Cal H}={\Cal
H}(G(\Q_{p^{n'}})// K_{p^{n'}})\leqno(10.6)$$ be the Hecke algebra
corresponding to the parahoric subgroup $K_{p^{n'}}$. It may be
conjectured that there exists an element $\phi_p^{n}\in {\Cal H}$ with the
following property. Let $g'_p\in G(\Q_{p^{n'}})$ be a representative of
$x'_p$. Then $${\roman{Contr}}_n(x)= \phi_p^n (g_p^{\prime -1}\delta
\sigma (g'_p))\ \ .\leqno(10.7)$$ Appealing to [K1], 1.5, we therefore
obtain the following group-theoretic expression for the contribution of
the admissible homomorphism $\varphi: {\Qfr} \to{\Cal G}_{\G}$ to the Lefschetz
fixed point formula over $\kappa_E^n$, $${\roman{Contr}}_n(\varphi)=v\cdot
O_h(\phi^p)\cdot TO_\delta(\phi_p^n)\ \ .\leqno(10.8)$$ Here $O_h(\phi^p)$
is the orbital integral over $h\in \G(\A_f^p)$ of the characteristic
function of $K^p$ and $TO_\delta(\phi_p^n)$ the twisted orbital integral
of $\phi_p^n$ over the twisted conjugacy class of $\delta\in
G(\Q_{p^{n'}})$. Furthermore, $v$ is a certain volume factor.

\medskip\noindent
For the function $\phi_p^n$ there is the following conjecture.

\medskip\noindent
{\bf Conjecture 10.3.} (Kottwitz) {\it Assume that $G$ splits over
$\Q_{p^{n'}}$. Let $K^{\circ}_{p^{n'}}$ be an Iwahori subgroup of
$G(\Q_{p^{n'}})$ contained in $K_{p^{n'}}$. Then $\phi_p^n$ is the image
of $p^{n'\cdot\langle\varrho,\mu\rangle}\cdot z_\mu$ under the homomorphism
of Hecke algebras $${\Cal
H}(G(\Q_{p^{n'}})//K^\circ_{p^{n'}})\longrightarrow {\Cal
H}(G(\Q_{p^{n'}})// K_{p^{n'}})\ \ .$$ Here $z_\mu$ denotes the Bernstein
function in the center of the Iwahori Hecke algebra associated to $\mu$,
comp.}\ [H2], 2.3.

\medskip\noindent
Recall that the center of ${\Cal H}(G(\Q_{p^{n'}})// K^\circ_{p^{n'}})$
has a basis as a $\C$-vector space formed by the Bernstein functions
$z_\lambda$, where $\lambda$ runs through the conjugacy classes of
one-parameter subgroups of $G$.

\medskip\noindent
In this direction we have the following facts.

\medskip\noindent
{\bf Theorem 10.4.} {\it Conjecture 10.3 holds in the following cases.
\item{(i)} {\rm (Haines, Ngo [HN2]):} $G=GL_n$ or $G=GSp_{2n}$.
\item{(ii)} {\rm (Haines [H2]):} $G$ is an inner form of $GL_n$ and
$\{\mu\}\ni\omega_1$ (the Drinfeld case).}

\medskip\noindent
It would be interesting to extend the statement of Kottwitz' conjecture to
the general case. Once this is done (and the corresponding conjecture
proved!) it remains to calculate the sum over all equivalence classes of
admissible homomorphisms $\varphi$ of the expressions (10.8) for
${\roman{Contr}}_n(\varphi)$. More precisely, one would like to replace
the twisted orbital integrals in (10.8) by an ordinary orbital integral of
a suitable function on $G(\Q_p)$ and compare the resulting expression with
the trace of a suitable function on $G({\bold A})$ in the automorphic
spectrum. When the Shimura variety is not projective, one also has to deal
with the contribution of the points on the boundary. Even when the Shimura
variety is projective, the phenomenon of $L$-indistinguishability
complicates the picture. But, at least these complications are of a
different nature from the ones addressed in this report. They are of a
group-theoretic nature, not of a geometric nature.

\bigskip

\subheading{Bibliography}

\bigskip\noindent

\ref{[BD]} Beilinson, A., Drinfeld, V.: Quantization of Hitchin's
integrable system and Hecke eigensheaves. Available at
http://www.math.uchicago.edu/~benzvi/

\ref{[BT1]} Bruhat, F., Tits, J.: Groupes r\'eductifs sur un corps local. Inst. Hautes
Etudes Sci. Publ. Math.  {\bf 41} (1972), 5--251.

\ref{[BT2]} Bruhat, F., Tits, J.: Groupes r\'eductifs sur un corps local. II.
Sch\'emas en groupes. Existence d'une donn\'ee radicielle valu\'ee.  Inst.
Hautes Etudes Sci. Publ. Math. {\bf 60} (1984), 197--376.

\ref{[C1]} Chai, C.-L.: Every ordinary symplectic isogeny class in positive
characteristic is dense in the moduli. Invent. Math. {\bf 121} (1995),
439--479.

\ref{[C2]} Chai, C.-L.: Newton polygons as lattice points. Amer. J. Math. {\bf 122} (2000),
967--990.

\ref{[CN1]} Chai, C.-L., Norman, P.: Bad reduction of the Siegel moduli scheme of
genus two with $\Gamma\sb 0(p)$-level structure. Amer. J. Math. {\bf 112}
(1990), 1003--1071.

\ref{[CN2]} Chai, C.-L., Norman, P.: Singularities of the $\Gamma\sb 0(p)$-level
structure. J. Algebraic Geom. {\bf 1} (1992),  251--278.

\ref{[DL]} Deligne, P., Lusztig, G.: Representations of reductive groups over finite
fields. Ann. of Math. {\bf 103} (1976), 103--161.

\ref{[DR]} Deligne, P., Rapoport, M.: Les sch\'emas de modules de courbes
elliptiques. Modular functions of one variable, II (Proc. Internat. Summer
School, Univ. Antwerp, Antwerp, 1972), pp. 143--316. Lecture Notes in
Math. {\bf 349}, Springer, Berlin, 1973.

\ref{[DM]} Digne, F., Michel, J.: Representations of finite groups of Lie type.
London Mathematical Society Student Texts, {\bf 21}. Cambridge University Press,
Cambridge, 1991. iv+159 pp.

\ref{[D]} Drinfeld, V. G.: Coverings of $p$-adic symmetric domains. (Russian)
Funk\-cio\-nal. Anal. i Prilo\v zen. {\bf 10} (1976),  29--40.

\ref{[F1]} Faltings, G.: Explicit resolution of local singularities of moduli-spaces.
J. Reine Angew. Math. {\bf 483} (1997), 183--196.

\ref{[F2]} Faltings, G.: Toroidal resolutions for some matrix singularities. Moduli
of abelian varieties (Texel Island, 1999), 157--184, Progr. Math.\ {\bf 195},
Birkh\"auser, Basel, 2001.

\ref{[Ge]} Genestier, A.: Un mod\`ele semi-stable de la vari\'et\'e de Siegel de
genre 3 avec structures de niveau de type $\Gamma\sb 0(p)$. Compositio
Math. {\bf 123} (2000),  303--328.

\ref{[G1]} G\"ortz, U.: On the flatness of models of certain Shimura varieties of
PEL-type. Math. Ann. {\bf 321} (2001), 689-727.

\ref{[G2]} G\"ortz, U.:  On the flatness of local models for the symplectic group. To
appear in Adv.\ Math., math.AG/0011202

\ref{[G3]} G\"ortz, U.: Topological flatness of local models in the ramified
case. In preparation.


\ref{[H1]} Haines, Th.: On matrix coefficients of the Satake isomorphism: complements
to the paper of M. Rapoport. Manuscripta Math. {\bf 101} (2000),  167--174.

\ref{[H2]} Haines, Th.: Test functions for Shimura varieties: the Drinfeld case. Duke
Math. J. {\bf 106} (2001), 19--40.

\ref{[H3]} Haines, Th.: The combinatorics of Bernstein functions. Trans. Amer. Math.
Soc. {\bf 353} (2001),  1251--1278.

\ref{[HN1]} Haines, Th., Ngo, B.C.: Alcoves associated to special fibers of local
models. To appear in American J.\ Math.,
 math.RT/0103048

\ref{[HN2]} Haines, Th., Ngo, B.C.: Nearby cycles for local models of some Shimura
varieties. To appear in Compositio Math.,
 math.AG/0103047

\ref{[HR]} Haines, Th., Rapoport, M.: Paper in preparation.

\ref{[IM]} Iwahori, N., Matsumoto, H.: On some Bruhat decomposition and the structure
of the Hecke rings of $p$-adic Chevalley groups. Inst. Hautes Etudes Sci.
Publ. Math. {\bf 25} (1965), 5--48.

\ref{[J]} de Jong, A.: The moduli spaces of principally polarized abelian varieties
with $\Gamma\sb 0(p)$-level structure. J. Algebraic Geom. {\bf 2} (1993),
667--688.

\ref{[JO]} de Jong, A., Oort, F.: Purity of the stratification by Newton polygons. J.
Amer. Math. Soc. {\bf 13} (2000), 209--241.

\ref{[Ka]} Kaiser, Ch.: Ein getwistetes fundamentales Lemma f\"ur die
${\roman{GSp}}\sb 4$. Bonner Mathematische Schriften, {\bf 303}. Universit\"at
Bonn, Mathematisches Institut, Bonn, 1997. 71 pp.

\ref{[K1]} Kottwitz, R.: Shimura varieties and twisted orbital integrals. Math. Ann.
{\bf 269} (1984), 287--300.

\ref{[K2]} Kottwitz, R.: Isocrystals with additional structure. Compositio Math. {\bf
56}
(1985), 201--220.

\ref{[K3]} Kottwitz, R.: Points on some Shimura varieties over finite fields. J.
Amer. Math. Soc. {\bf 5} (1992), 373--444.

\ref{[K4]} Kottwitz, R.: Isocrystals with additional structure. II. Compositio Math.
{\bf 109} (1997), 255--339.

\ref{[KR1]} Kottwitz, R., Rapoport, M.: Minuscule alcoves for ${\roman{GL}}\sb n$ and
${\roman{GSp}}\sb {2n}$. Manu\-scrip\-ta Math. {\bf 102} (2000), 403--428.

\ref{[KR2]} Kottwitz, R., Rapoport, M.: On the existence of $F$-crystals.
\hfill\break
math.NT/0202229

\ref{[L]} Lafforgue, L.: Pavages de poly\`edres, recollement de cellules
de Schubert et compactification d'espaces de configurations. In
preparation.

\ref{[LR]} Langlands, R., Rapoport, M.: Shimuravariet\"aten und Gerben. J. Reine
Angew. Math. {\bf 378} (1987), 113--220.

\ref{[LO]} Li, K.-Z., Oort, F.: Moduli of supersingular abelian varieties. Lecture
Notes in Mathematics {\bf 1680}. Springer-Verlag, Berlin, 1998. iv+116 pp.

\ref{[M1]} Milne, J.: The conjecture of Langlands and Rapoport for Siegel modular
varieties. Bull. Amer. Math. Soc. (N.S.) {\bf 24} (1991), 335--341.

\ref{[M2]} Milne, J.: The points on a Shimura variety modulo a prime of good
reduction. The zeta functions of Picard modular surfaces, 151--253, Univ.
Montr\'eal, Montreal, QC, 1992.

\ref{[M3]} Milne, J.: Motives over finite fields. Motives (Seattle, WA, 1991),
401--459, Proc. Sympos. Pure Math. {\bf 55}, Part 1, Amer. Math. Soc.,
Providence, RI, 1994.

\ref{[NG]} Ngo, B.C., Genestier, A.: Alcoves et p-rang des vari\'et\'es ab\'eliennes.
\hfill\break
math.AG/0107223

\ref{[O]} Oort, F.: Newton polygon strata in the moduli space of abelian varieties.
Moduli of abelian varieties (Texel Island, 1999), 417--440, Progr. Math.,
{\bf 195}, Birkh\"auser, Basel, 2001.

\ref{[P]} Pappas, G.: On the arithmetic moduli schemes of PEL Shimura varieties. J.\
Alg.\ Geom.\ {\bf 9} (2000), 577--605.

\ref{[PR1]} Pappas, G., Rapoport, M.: Local models in the ramified case I. The
EL-case. To appear in J.\ Alg.\ Geom., math.AG/0006222

\ref{[PR2]} Pappas, G., Rapoport, M.: Local models in the ramified case II. Splitting
models. In preparation.

\ref{[Pf]} Pfau, M.: The conjecture of Langlands and Rapoport for certain Shimura
varieties of non-rational weight. J. Reine Angew. Math. {\bf 471} (1996),
165--199.

\ref{[R1]} Rapoport, M.: On the local zeta function of quaternionic Shimura varieties
with bad reduction. Math. Ann. {\bf 279} (1988), 673--697.

\ref{[R2]} Rapoport, M.: On the bad reduction of Shimura varieties. Automorphic
forms, Shimura varieties, and $L$-functions, Vol. II (Ann Arbor, MI,
1988), 253--321, Perspect. Math.\ {\bf 11}, Academic Press, Boston, MA, 1990.

\ref{[R3]} Rapoport, M.: A positivity property of the Satake isomorphism. Manu\-scrip\-ta
Math. {\bf 101} (2000), 153--166.

\ref{[RR]} Rapoport, M., Richartz, M.: On the classification and specialization of
$F$-isocrystals with additional structure. Compositio Math. {\bf 103} (1996),
153--181.

\ref{[RZ1]} Rapoport, M., Zink, Th.: \"Uber die lokale Zetafunktion von
Shimuravariet\"aten. Monodromiefiltration und verschwindende Zyklen in
ungleicher Charakteristik. Invent. Math. {\bf 68} (1982), 21--101.

\ref{[RZ2]} Rapoport, M., Zink, Th.: Period spaces for $p$-divisible groups. Annals of
Mathematics Studies, {\bf 141}. Princeton University Press, Princeton, NJ, 1996.
xxii+324 pp.

\ref{[Re1]} Reimann, H.: The semi-simple zeta function of quaternionic Shimura
varieties. Lecture Notes in Mathematics, {\bf 1657}. Springer-Verlag, Berlin,
1997. viii+143 pp.

\ref{[Re2]} Reimann, H.: On the zeta function of quaternionic Shimura
varieties. Math.\ Ann.\ {\bf 317} (2000), 41--55.

\ref{[Re3]} Reimann, H.: Reduction of Shimura varieties at parahoric
levels. Manu\-scrip\-ta Math.\ {\bf 107} (2002), 355--390.
\ref{[Ri]} Richartz, M.: Klassifikation von selbstdualen
Dieudonn\'egittern in einem dreidimensionalen polarisierten
supersingul\"aren Isokristall. Diss. Universit\"at Bonn, 1998.

\ref{[S]} Stamm, H.: On the reduction of the Hilbert-Blumenthal-moduli scheme with
$\Gamma\sb 0(p)$-level structure. Forum Math. {\bf 9} (1997),  405--455.

\ref{[T]} Tits, J.: Reductive groups over local fields.  Proc. Sympos. Pure Math. {\bf 33}, Part 1, pp. 29--69,
 Amer. Math. Soc.,
Providence, R.I., 1979.

\ref{[W]} Wedhorn, T.: Ordinariness in good reductions of Shimura varieties of
PEL-type. Ann. Sci. Ecole Norm. Sup. {\bf 32} (1999),  575--618.

\end